\documentclass[11pt]{article}

\usepackage{iftex}
\ifPDFTeX
  \usepackage[utf8]{inputenc}
\fi
\usepackage[T1]{fontenc}
\usepackage{lmodern}
\usepackage[a4paper,margin=2.5cm]{geometry}
\usepackage{amsmath,amssymb,amsthm,bm}
\usepackage{stmaryrd}
\usepackage{empheq}
\usepackage{enumitem}
\usepackage{graphicx}
\usepackage{subcaption}
\usepackage{ifpdf}
\usepackage{hyperref}
\usepackage[capitalise,noabbrev]{cleveref}

\numberwithin{equation}{section}

\theoremstyle{plain}
\newtheorem{theorem}{Theorem}[section]
\newtheorem{lemma}[theorem]{Lemma}
\newtheorem{proposition}[theorem]{Proposition}
\newtheorem{corollary}[theorem]{Corollary}
\theoremstyle{definition}
\newtheorem{definition}[theorem]{Definition}

\newcommand{\email}[1]{\href{mailto:#1}{\nolinkurl{#1}}}
\newenvironment{keywords}{\par\medskip\noindent\small\textbf{Keywords.}\ }{\par}
\newenvironment{MSCcodes}{\par\noindent\small\textbf{MSC codes.}\ }{\par}

\newcommand{\R}{\mathbb{R}}
\newcommand{\Z}{\mathbb{Z}}
\newcommand{\N}{\mathbb{N}}
\newcommand{\e}{\varepsilon}

\newcommand{\di}[1]{\,\mathrm{d}#1}

\newcommand{\dive}{\operatorname{div}}

\newcommand{\jump}[1]{\llbracket #1\rrbracket}

\newcommand{\indexI}{\mathcal{I}_\e}

\newcommand{\bh}{\bm{h}}

\setlist[itemize]{leftmargin=2em}
\setlist[enumerate]{leftmargin=*,labelindent=\parindent}
\setlist[enumerate]{leftmargin=2em}

\newcommand{\systemlabel}[2]{\hypertarget{#1}{\ensuremath{(#2)}}}
\newcommand{\systemref}[2]{\hyperlink{#1}{\ensuremath{(#2)}}}
\DeclareMathOperator{\dist}{dist}

\theoremstyle{remark}
\newtheorem{remark}[theorem]{Remark}

\crefname{hypothesis}{Hypothesis}{Hypotheses}

\crefname{example}{Example}{Examples}
\allowdisplaybreaks[1]
\title{Homogenization with Shape-Evolving Microstructures:\\
A Two-Phase Membrane Problem}

\author{Michael Eden\thanks{Faculty of Mathematics, University of Regensburg, Germany
  (\email{michael.eden@ur.de}).}}

\ifpdf
\hypersetup{
  pdftitle={Homogenization with Shape-Evolving Microstructures: A Two-Phase Membrane Problem},
  pdfauthor={M. Eden}
}
\fi

\begin{document}

\maketitle

\begin{abstract}
We study a diffusion problem in a highly heterogeneous periodic medium (characterized by a small scale parameter $0<\e\ll1$) where small inclusions are separated from the surrounding connected area by solute-impermeable membranes.
These membranes are allowed to grow and deform under surface tension and osmotic pressure (driven by the concentration difference across the membrane).
We prove existence of weak solutions on a time interval which is independent of the microscopic scale via a Schauder fixed-point argument by establishing spatial boundary Hölder estimates.
For the homogenization limit $\e\to0$, we establish strong compactness of the outer concentration via time translates and then derive a stability estimate for the driven inner cell problems.
The limit couples a macroscopic conservation law, with shape-dependent effective diffusion, to inner diffusion and membrane motion in each cell.
\end{abstract}



\section{Introduction}\label{sec:introduction}

Transport of chemical substances through heterogeneous media can heavily depend on the geometry of the microstructure and when this geometry changes in response to the transported substances, the effective material properties may change as well.
One natural example of such a feedback is osmotic swelling:
local concentration differences across a \textit{semipermeable membrane}\footnote{Here, this means permeable to the solvent but impermeable to the dissolved substances.} exert an osmotic pressure that drives the membrane's motion, while surface tension opposes the increase of interfacial area.
Free-boundary models describing this mechanism have been studied for one- and two-phase configurations in \cite{LippothProkert2012,LippothProkert2014};
for corresponding models using phase-field approximations, we refer to \cite{Raetz2016,TangLiLowengrubWise2023}.
In a medium containing many small inclusions, each separated from a connected surrounding phase via such a membrane, an effective description of the transport processes should take into account the size and shape of the membranes without resolving every individual membrane at the macroscopic scale.
Such heterogeneous configurations can be found, e.g., in biological tissues, where cell volume changes in response to osmotic gradients, in suspensions of vesicles, or in swelling hydrogels \cite{HoffmannLambertPedersen2009,TanakaFillmore1979, TangLiLowengrubWise2023}.

In this work, we investigate a two-phase diffusion problem in a connected outer region containing periodically distributed inclusions.
The separating membranes of these inclusions are impermeable to the solute, and their normal velocity is determined by the local concentration difference and the mean curvature.
We use a reduced osmosis model where we neglect the convective transport within the phases; this is a common simplification in osmosis models when diffusion dominates convection on the scale of the inclusions \cite{Zaal2012,LippothPeletierProkert2016}.

\paragraph{Microscopic model and effective equations.}
Let $\Omega_\e^{+}(t)$ denote the connected outer phase, $\Omega_\e^{-}(t)$ the union of inclusions, and $\Gamma_\e(t)$ their common interface (with the scale parameter $0<\e\ll1$ denoting the ratio of the microscopic to the macroscopic length scale).
The normal $n_\e$ points out of the inclusions, and $v_\e$ and $H_\e$ denote the normal velocity and mean curvature with this orientation.
The essential equations of the microscopic problem \systemref{sys:osmosis}{P_\e} are
\[
\begin{aligned}
 \partial_tu_\e^{\pm}-d_\e^{\pm}\Delta u_\e^{\pm}&=f_\e^{\pm}
     &&\text{in }\Omega_\e^{\pm}(t),\\
 d_\e^{\pm}\nabla u_\e^{\pm}\cdot n_\e+u_\e^{\pm}v_\e&=0
     &&\text{on }\Gamma_\e(t),\\
 v_\e&=-\e^2\sigma H_\e+\e\lambda(u_\e^{-}-u_\e^{+})
     &&\text{on }\Gamma_\e(t),
\end{aligned}
\]
where $d_\e^{+}=d^{+}$, $d_\e^{-}=\e^2d^{-}$, and $d^{+},d^{-},\sigma,\lambda>0$ are fixed.
This $\e$-scaling is the classical double-porosity scaling \cite{ArbogastDouglasHornung1990,HornungJaegerMikelic1994}, also used for inclusions of macroscopically varying size in \cite{vanNoordenMuntean2011,MunteanvanNoorden2013}.
The zero flux condition encodes the impermeability of the membranes and the system is completed by initial data and a no-flux condition at the fixed outer boundary (see \systemref{sys:osmosis}{P_\e}).

In the limit $\e\to0$, the outer concentration $u^{+}(t,x)$ depends only on the macroscopic position, whereas the inner concentration $u^{-}(t,x,y)$ is defined on the moving inclusion $Y^{-}(t,x)$ inside the unit cell $Y=(0,1)^3$.
We write $Y^{+}(t,x)=Y\setminus\overline{Y^{-}(t,x)}$ and $\Gamma(t,x)=\partial Y^{-}(t,x)$.
With $V$ and $H$ again measured using the outward inclusion normal, the corresponding equations of the limit problem \systemref{sys:homogenized-moving}{P_0} are
\[
\begin{aligned}
 \partial_t\bigl(\theta^{+}(h)u^{+}\bigr)
       -\dive_x\bigl(D(h)\nabla_xu^{+}\bigr)&=\theta^{+}(h)f^{+}
       &&\text{in }\Omega,\\
 \partial_tu^{-}-d^{-}\Delta_y u^{-}&=f^{-}
       &&\text{in }Y^{-}(t,x),\\
 d^{-}\nabla_y u^{-}\cdot n+u^{-}V&=0
       &&\text{on }\Gamma(t,x),\\
 V&=-\sigma H+\lambda(u^{-}-u^{+})
       &&\text{on }\Gamma(t,x).
\end{aligned}
\]
Here, $h=h(t,x,\gamma)$ describes the cell interface as a normal graph over the reference membrane.
The outer volume fraction is $\theta^{+}(h)=|Y^{+}(t,x)|$ and $D(h)$ is the effective diffusion tensor obtained via a typical elliptic cell problem over $Y^{+}(t,x)$.

Our main results justify this passage on a common local time interval:
\Cref{thm:prepared_schauder} provides weak microscopic solutions on a time horizon $I=(0,T)$ (where $T>0$ is independent of $\e$) with uniformly bounded concentrations and uniformly controlled cell shapes and velocities.
Then, \cref{thm:homogenization} proves subsequence convergence of these solutions to a solution of the coupled limit above.
After the reference-domain transformations and extensions, the outer concentration converges strongly in every finite space-time $L^q$.
The unfolded inner concentrations and heights converge strongly as well, the latter in $L^q(\Omega;C(\overline{I};C^{2+\beta'}(\Gamma)))$ for $1\le q<\infty$ and some $\beta'>0$.
The precise data assumptions, convergence statements, and weak interpretation of the moving-cell equations are given in \cref{sec:setting}.

\paragraph{Previous results and relation to the literature.}
A common approach to homogenization with evolving microstructures is to transform the moving domains to a fixed reference configuration \cite{Peter2009}.
This has been used for thermoelasticity with prescribed phase transformations \cite{EdenMuntean2017}, for motions generated by a prescribed normal velocity via a height function \cite{Eden2019}, for nonlinear transport on evolving domains \cite{GahnNeussRaduPop2021}, and for Stokes flow, where the change of the local pore volume produces an additional macroscopic source term \cite{WiedemannPeter2024}.
In all these works above, the geometry may vary substantially over time, but its evolution is given \textit{a priori} and does not depend on the transport process.

For a free microscopic geometry, the difficulty is to control the transport and the motion simultaneously.
A porous medium with swelling spherical granules, whose radii follow a moisture-dependent kinetic law, was analyzed in \cite{FasanoMikelic2002}.
A coupled reaction--diffusion model with evolving pore volume was analyzed and homogenized in \cite{Peter2009}.
In cell tissue, a feedback of the transported species on the microstructure through concentration-dependent elastic moduli leads to effective coefficients depending on the concentration \cite{JaegerMikelicNeussRadu2011}.
More recently, rigorous results for concentration-dependent growth and dissolution of spherical inclusions were obtained in \cite{GahnPop2023,WiedemannPeter2023}.
In a similar setup, Stokes flow together with reactive transport is considered in \cite{GahnPeterPopWiedemann2024}, again with spherical grains described by their radii.
In \cite{EdenMuntean2026}, the evolution of general smooth reference inclusions (not just spherical ones) is considered in a thermoelastic phase-transition model, but the height is still constant for each individual reference interface (thereby not allowing shape changes).
An alternative approach is given by diffuse-interface descriptions, in which a phase-field variable approximates the evolving pore geometry instead of tracking the
interface exactly.
For precipitation and dissolution in porous media, such models have been upscaled by
\textit{formal} asymptotic expansions (e.g., \cite{BringedalVonWolffPop2020,vonWolffPop2022}),
leading to cell problems for the phase field that are structurally similar to the cell
problems obtained in our problem.

In both the microscopic problem and the limit, the height function describing the interface motion depends on the surface points $\gamma\in\Gamma$, and both the concentration difference and the curvature in the velocity law may vary along the interface.
Therefore, the evolving geometry is not restricted to a family parametrized by a radius or a single cellwise height; this is the sense in which we consider genuine shape changes of small normal graphs over smooth reference interfaces.
To the best of our knowledge, this is the first rigorous homogenization result for a free microscopic geometry whose interfaces change their shape, not only their size, in response to the transported quantities.

\paragraph{Analytical strategy and main contributions.}
The two central difficulties are $(i)$ to find an existence interval independent of $\e$ and $(ii)$ to establish sufficiently strong convergence of the geometry to be able to pass to the limit.
\begin{itemize}
\item[$(i)$] We note first that local existence on some interval $(0,T_\e)$ (for fixed $\e$) follows from \cite{LippothProkert2014}, at least for sufficiently regular and compatible data.
However, $T_\e$ may degenerate as $\e\to0$, since the number of interfaces grows and their distances decrease.
In our approach, we derive uniform spatial Hölder bounds for the concentration traces on the membranes, which are enough to control the interface motion via maximal parabolic regularity \cite{DenkPruessZacher2008,PruessSimonett2016}.
A Schauder fixed-point argument then gives an $\e$-independent existence interval.
The construction is motivated by the fast-diffusion analysis of \cite{Lippoth2016}.
\item[$(ii)$] We first establish strong compactness of the outer concentration by a time-translation argument for the weighted mass $J_\e\widehat u_\e^{+}$, where $J_\e$ is the Jacobian determinant of the transformation.
This is a standard argument in homogenization with evolving microstructures, see, e.g., \cite{GahnNeussRaduPop2021,WiedemannPeter2023,GahnPeterPopWiedemann2024}.
For the inner cells, we prove a general stability estimate for the cell problem with respect to forcing data (see \cref{thm:cell_stability}).
Using this stability result, we are able to transfer the strong convergence of the outer concentration and the data to the inner concentrations and the heights, and hence to the geometry.
This is then enough to pass to the limit.
\end{itemize}

\paragraph{Limitations and outlook.}
Our results are local in time: 
the heights must correspond to small normal graphs over the reference interface, so that the Hanzawa transformation stays a diffeomorphism and the maximal-regularity argument of \cref{lem:surface_spatial} applies.
Large deformations, contact between inclusions or with the cell boundary, and topological changes are therefore inherently excluded.
In addition, the homogenization requires strongly convergent data \ref{ass:macro_initial}--\ref{ass:macro_sources} to pass to the limit.
Uniqueness of the limit problem also remains open, although we expect weak–strong uniqueness to hold for sufficiently regular solutions; see \cref{rem:homogenization}.

On the other hand, our strategy does not rely on the specific form of the osmosis model beyond two structural features: surface tension makes the interface evolution parabolic, and the one-phase cell problem is stable with respect to the driving outer concentration.
We therefore expect it to extend, for instance, to models with additional bulk flow
\cite{LippothPeletierProkert2016}, to Stefan problems with Gibbs--Thomson correction and kinetic undercooling, to precipitation and dissolution processes in which the growth rate varies along the grain surface, cf.~\cite{GahnPop2023,GahnPeterPopWiedemann2024}, or to thermoelastic phase transitions with spatially varying interface heights, cf.~\cite{EdenMuntean2026}.

\section{Notation, model, and main results}
\label{sec:setting}

The superscripts $+$ and $-$ refer to the connected outer phase and the inner inclusions, respectively.
A formula containing $\pm$ is understood separately for each phase, with the same sign chosen throughout that formula; $\sum_{\pm}$ denotes the sum of these two expressions.
For a real quantity $z$, we write $\operatorname{pos}(z)=\max\{z,0\}$ and $\operatorname{neg}(z)=\max\{-z,0\}$ to denote the positive and negative part.
For $p\in[1,\infty)$, $m\in\N_0$, and $0<\alpha<1$, we use the usual Lebesgue ($L^p$), Sobolev ($H^m=W^{m,2}$), and Hölder spaces ($C^{m+\alpha}$).
For parabolic Hölder spaces, we use the standard notation $C^{\beta,\alpha}$, with the first exponent referring to time and the second to space.
Constants denoted by $C$ may change between occurrences; in uniform estimates they are independent of $\e$ and the cell index.

\begin{figure}[htbp]
    \centering
    \includegraphics[width=\linewidth]{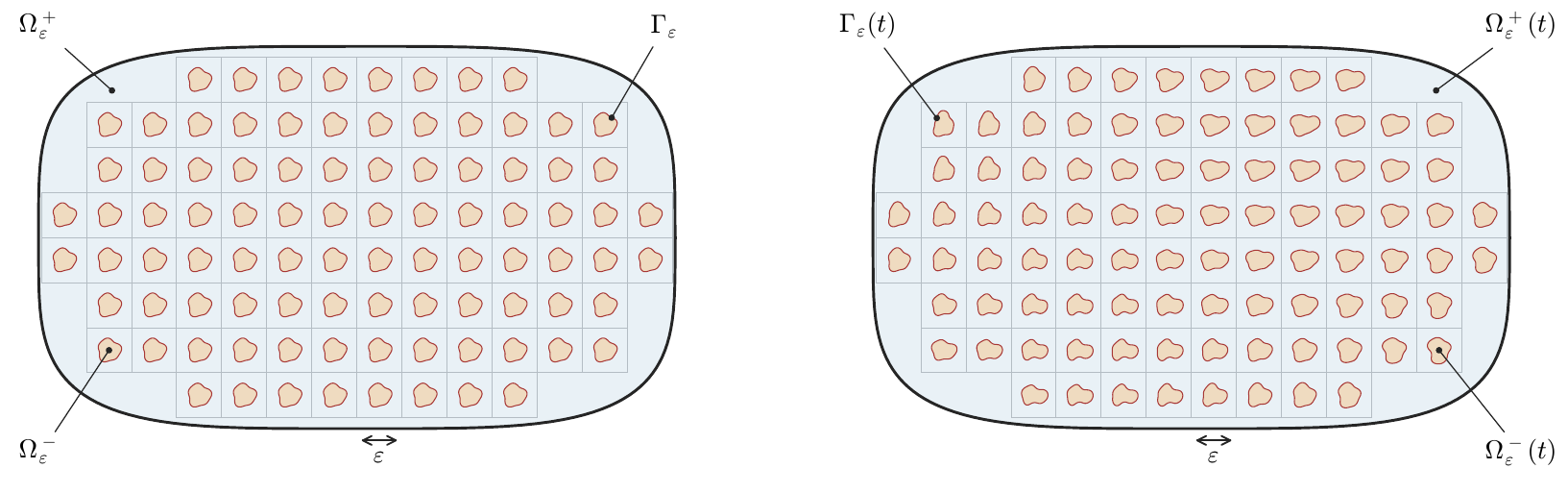}
    \caption{Left: initial configuration, in which all inclusions are translates of $\e Y^-$. 
    Right: configuration at a later time, in which each inclusion has deformed individually.}
    \label{fig:cell_variation}
\end{figure}

Let $\Omega\subset\R^3$ be a bounded Lipschitz domain and fix a data horizon $T^*>0$. 
We write $I^*:=(0,T^*)$ and $I_\tau:=(0,\tau)$ for $\tau\in(0, T^*)$. The symbol $T$ denotes the common existence horizon, independent of $\e$; once it has been chosen, we write $I:=I_T=(0,T)$. 
We use the reference cell configuration
\[
 Y=(0,1)^3,\qquad Y^{-}\Subset Y,\qquad
 Y^{+}=Y\setminus\overline{Y^{-}},\qquad \Gamma=\partial Y^{-},
\]
where $Y^{-}$ and $Y^{+}$ are connected and $\Gamma$ is $C^4$.\footnote{$C^2$ interfaces could in principle be handled by writing $\Gamma$ as a normal graph over a nearby smooth auxiliary surface \cite[Chapter~2]{PruessSimonett2016}.
Also, there is no principal obstacle to allowing \(Y^{-}\) to have several connected components.
We work in this more regular regime to simplify the notation and arguments.}
For a sequence $\e\downarrow0$, $\e\le1$, we set
\[
 \indexI=\{k\in\Z^3:\ \e(k+\overline{Y})\subset\Omega\}
\]
and note that $|\indexI|\simeq|\Omega|\e^{-3}$.
We set
\[
 \Omega_\e^{-}=\bigcup_{k\in \indexI}\e(k+Y^{-}),\qquad
 \Omega_\e^{+}=\Omega\setminus\overline{\Omega_\e^{-}},\qquad
 \Gamma_\e=\bigcup_{k\in \indexI}\e(k+\Gamma).
\]
We note that with this construction, $\partial\Omega_\e^-=\Gamma_\e$ and $\partial\Omega_\e^+=\Gamma_\e\cup\partial\Omega$.
Moreover, $\Omega_\e^-$ is also $C^4$-regular.
We denote the moving phases by $\Omega_\e^{\pm}(t)$.
In particular, $\Gamma_\e=\Gamma_\e(0)$ and $\Omega_\e^\pm=\Omega_\e^\pm(0)$.
We refer to \cref{fig:cell_variation} for a sketch of this geometric setup.

On the fixed reference interface $\Gamma$, we let $n_\Gamma$ denote the unit normal pointing out of $Y^{-}$ and, on the moving microscopic interface $\Gamma_\e(t)$, we let $n_\e=n_\e(t)$ denote the unit normal pointing out of $\Omega_\e^{-}(t)$.
We introduce the respective mean curvatures by
\[
 H_\Gamma=\dive_\Gamma n_\Gamma,
 \qquad
 H_\e(t)=\dive_{\Gamma_\e(t)}n_\e(t).
\]
with the convention that they are positive on spheres.
The scalar $v_\e$ is the normal velocity in the direction $n_\e$; positive velocity therefore increases the volume of the inner phase at the expense of the outer phase.
For quantities with traces on both sides of the membrane, our jump convention is $\jump{z}=z^{+}|_\Gamma-z^{-}|_\Gamma$, with the corresponding traces on $\Gamma_\e(t)$ in physical coordinates.

\subsection{Microscopic model and assumptions}\label{ssec:assumptions}

Let $u_\e^{\pm}$ denote the concentration of the dissolved species in the respective phase.
We fix parameters $d^{+},d^{-},\sigma,\lambda>0$ and set
\[
 d_\e^{+}=d^{+},\qquad d_\e^{-}=\e^2d^{-}.
\]
Then, for prescribed sources $f_\e^{\pm}$ and initial values $u_{\e,0}^\pm$, the microscopic two-phase membrane problem is given by
\begin{empheq}[left={\systemlabel{sys:osmosis}{P_\e}\quad\empheqlbrace}]{alignat=2}
 \partial_tu_\e^{\pm}-d_\e^{\pm}\Delta u_\e^{\pm}&=f_\e^{\pm}
 &&\quad\text{in }\Omega_\e^{\pm}(t),
 \label{osmosis_problem:bulk}\\
 d_\e^{\pm}\nabla u_\e^{\pm}\cdot n_\e+u_\e^{\pm}v_\e&=0
 &&\quad\text{on }\Gamma_\e(t),\label{osmosis_problem:flux:plus}\\
 v_\e&=-\e^2\sigma H_\e-\e\lambda\jump{u_\e}
 &&\quad\text{on }\Gamma_\e(t),\label{osmosis_problem:velocity}\\
 \nabla u_\e^{+}\cdot\nu&=0
 &&\quad\text{on }\partial\Omega,\label{osmosis_problem:outer}\\
 u_\e^{\pm}(0)&=u_{\e,0}^{\pm},\qquad \Gamma_\e(0)=\Gamma_\e.
 &&\label{osmosis_problem:initial}
\end{empheq}
Under the transformation $x=\e(k+y)$, the curvature and velocity scale as
\begin{equation}\label{eq:scaling}
 H_\e=\e^{-1}H,\qquad v_\e=\e V,\qquad
 \nabla_x=\e^{-1}\nabla_y.
\end{equation}
Here $H$ and $V$ are the physical curvature and normal velocity on the rescaled moving membrane, measured with respect to its outward normal.

Throughout the analysis we impose the following assumptions on the geometry and data.
The limit assumptions \ref{ass:macro_initial} and \ref{ass:macro_sources} are used only in the homogenization argument.
\begin{enumerate}[label=\textbf{(A\arabic*)},ref=(A\arabic*)]
\item\label{ass:geometry} \emph{(Initial geometry.)} The reference and outer geometries satisfy the conditions above.
The constants $d^{\pm},\sigma,\lambda>0$ are strictly positive.
\item\label{ass:initial_data} \emph{(Initial data.)} There are $u^*,U^*>0$ such that, for every $\e>0$,
\[
 u_{\e,0}^\pm\in L^\infty(\Omega_\e^\pm),
 \qquad 0\le u_{\e,0}^\pm\le u^*,
 \qquad
 \sup_{k\in\indexI}
  \bigl\|u_{\e,0}^\pm\bigl(\e(k+\cdot)\bigr)\bigr\|_{W^{2,\infty}(Y^\pm)}\le U^* .
\]
\item\label{ass:sources} \emph{(Ambient sources.)} The sources are defined on the whole domain, $f_\e^\pm\in L^\infty(I^*\times\Omega)$, and there is $F^*>0$ such that
\begin{equation*}
 \sup_{\e>0}\Bigl(\|f_\e^+\|_{L^\infty(I^*;W^{1,\infty}(\Omega))}+
 \|f_\e^-\|_{L^\infty(I^*\times\Omega)}
 +\sup_{k\in\indexI}\bigl\|f_\e^-\bigl(\cdot,\e(k+\cdot)\bigr)\bigr\|_{L^\infty(I^*;W^{1,\infty}(Y))}\Bigr)\le F^* .
\end{equation*}
\end{enumerate}
\begin{enumerate}[label=\textbf{(H\arabic*)},ref=(H\arabic*)]
\item\label{ass:macro_initial} \emph{(Macroscopic initial data.)}
There are $u_0^+\in L^2(\Omega)$ and $u_0^-\in L^2(\Omega\times Y^-)$ such that
\begin{gather*}
 P_\e^+u_{\e,0}^+\to u_0^+
 \ \text{ strongly in }L^2(\Omega),\quad
 \mathcal T_\e u_{\e,0}^-\to u_0^-
 \ \text{ strongly in }L^2(\Omega\times Y^{-}).
\end{gather*}
\item\label{ass:macro_sources} \emph{(Macroscopic sources.)}
There are $f^+\in L^2(I^*\times\Omega)$ and $f^-\in L^2(I^*\times\Omega\times Y)$ such that
\begin{gather*}
 f_\e^+\to f^+ \ \text{ strong in }L^2(I^*\times\Omega),\ 
 \mathcal T_\e f_\e^-\to f^- \ \text{ strongly in }L^2(I^*\times\Omega\times Y).
\end{gather*}
\end{enumerate}
Here $P_\e^+$ is the extension operator of \cref{lem:hom_extension} and $\mathcal T_\e$ is the periodic unfolding operator of Appendix~\ref{app:unfolding}.
Since \ref{ass:initial_data} controls $u_{\e,0}^+$ only cellwise, it may jump across common cell faces and need not belong to $H^1(\Omega_\e^+)$.
By the $L^2$ estimate in \cref{lem:hom_extension} and density, $P_\e^+$ extends uniquely to a bounded operator on $L^2(\Omega_\e^+)$, and \ref{ass:macro_initial} refers to this extension.

\begin{remark}\label[remark]{rem:assumptions}
The limits in \ref{ass:macro_initial}--\ref{ass:macro_sources} inherit the microscopic
 bounds: $0\le u_0^\pm\le u^*$, $\|u_0^-(x,\cdot)\|_{W^{2,\infty}(Y^-)}\le U^*$ for
 a.e.\ $x$, $\|f^+\|_{L^\infty(I^*;W^{1,\infty}(\Omega))}\le F^*$, and $f^-\in L^\infty(I^*\times\Omega;W^{1,\infty}(Y))$.
 By interpolation, the convergences of $\mathcal T_\e u_{\e,0}^-$ and $\mathcal T_\e f_\e^-$
 improve to $L^q(\Omega;C(\overline{Y^-}))$ and $L^q(I^*\times\Omega;C(\overline{Y}))$ for every
 $q<\infty$. Assumption \ref{ass:macro_sources} holds for $\e$-independent sources and for
 cell-scale oscillations $f_\e^-(t,x)=g(t,x,x/\e)$ with $g$ bounded, Lipschitz in $x$, and Lipschitz and $Y$-periodic in its
 last argument, with $f^-=g$.
\end{remark}

\subsection{Fixed-domain formulation and solution concept}
\label{ssec:fixed_domain}

To be able to formulate the problem on fixed reference domains, we parametrize the moving membranes by normal graphs; often called \textit{Hanzawa} transformation~\cite{Hanzawa1981}.
To that end, let $U_\Gamma\Subset Y$ be a tubular neighborhood of thickness $a_\Gamma>0$ of $\Gamma$,\footnote{We could allow different inner and outer thicknesses, but this would only complicate the notation. 
The thickness $a_\Gamma$ is determined by the geometry of \(\Gamma\) and its position inside \(Y\).
See, e.g., \cite[Section~2.3]{PruessSimonett2016} for tubular neighborhoods and normal coordinates of smooth hypersurfaces.} with coordinates
\[
 \Lambda(\gamma,\zeta)=\gamma+\zeta n_\Gamma(\gamma),\qquad |\zeta|<a_\Gamma,\ \gamma\in\Gamma.
\]
Choose $\chi\in C_c^\infty(-a_\Gamma,a_\Gamma)$ with $\chi\equiv1$ on $(-\frac{a_\Gamma}2,\frac{a_\Gamma}2)$.
For $h\colon\Gamma\to\R$, we set
\begin{equation}\label{eq:hanzawa_cell}
\Theta_h\colon Y\to Y,
    \qquad
    \Theta_h(y)
    =
    \begin{cases}
        \Lambda(\gamma,\zeta)+\chi(\zeta)h(\gamma)n_\Gamma(\gamma),
        & y=\Lambda(\gamma,\zeta)\in U_\Gamma,\\[1mm]
        y,
        & y\in Y\setminus U_\Gamma
    \end{cases}
\end{equation}
and choose a fixed number $\eta^*>0$ sufficiently small such that, for every $h\in C^2(\Gamma)$ with $\|h\|_{C^1(\Gamma)}\le\eta^*$, the map $\Theta_h$ is a diffeomorphism of $\overline Y$; we call such heights \emph{admissible}.
We write
\[
 \Gamma[h]=\Theta_h(\Gamma)=\{\gamma+h(\gamma)n_\Gamma(\gamma)\ : \ \gamma\in\Gamma\},\qquad
 Y^{\pm}[h]=\Theta_h(Y^{\pm}).
\]
For admissible $h$, $\Gamma[h]$ is a compact embedded $C^2$-surface
\begin{figure}[htbp]
    \centering
    \includegraphics[width=0.7\linewidth]{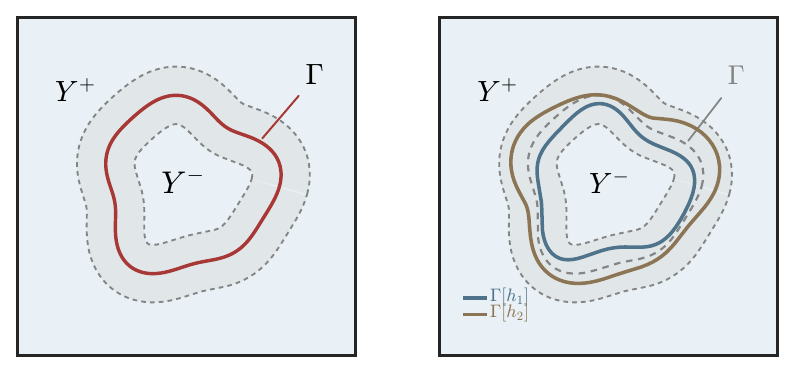}
    \caption{$(a)$ Reference cell $Y$ with tubular neighborhood, $(b)$ Cell interfaces for two different height functions $h_1,h_2$}
    \label{fig:cell_geometries}
\end{figure}
and we introduce the cell coefficients
\begin{equation}\label{eq:geometric_coefficients}
 J_h=\det D\Theta_h,\qquad
 A_h=J_hD\Theta_h^{-1}D\Theta_h^{-T},\qquad
 N_hz=\frac{D\Theta_h^{-1}D\Theta_h^{-T}\nabla z\cdot n_\Gamma}
 {|D\Theta_h^{-T}n_\Gamma|},
\end{equation}
where $D\Theta_h^{-1}=(D\Theta_h)^{-1}$ and $D\Theta_h^{-T}$ is its transpose.
Let $X_h(\gamma)=\gamma+h(\gamma)n_\Gamma(\gamma)$ denote the restriction of $\Theta_h$ to $\Gamma$,
let $n[h]$ be the outward unit normal on $\Gamma[h]$, and define its reference pullback by
\[
 \widehat n_h:=n[h]\circ X_h
 =\frac{D\Theta_h^{-T}n_\Gamma}{|D\Theta_h^{-T}n_\Gamma|},
 \qquad \beta(h):=n_\Gamma\cdot\widehat n_h.
\]
The factor $\beta(h)$ measures the tilt of the normal due to the deformation.
It equals one for constant $h$ and is bounded below by a fixed $\beta_*>0$ for all admissible $h$.
The physical cell curvature and its reference pullback are
\[
 H[h]:=\dive_{\Gamma[h]}n[h],\qquad
 \widehat H[h]:=H[h]\circ X_h,
\]
so that $\widehat H[h]$ and $\widehat n_h$ live on the fixed interface $\Gamma$, and
$\widehat H[0]=H_\Gamma$, $\widehat n_0=n_\Gamma$.
For the surface Jacobian, we have
\begin{equation}\label{eq:surface_jacobian}
 m_h=|D\Theta_h^{-T}n_\Gamma|,\qquad j_h=J_hm_h,
 \qquad \int_{\Gamma[h]}g\di S
       =\int_\Gamma(g\circ X_h)j_h\di S.
\end{equation}
We use subscripts for the geometric coefficients $J_h$, $A_h$, $j_h$, and square brackets for the nonlinear curvature operator $\widehat H[h]$.

Let now 
\begin{equation}\label{eq:defintion_XI}
h\in\mathcal X(I):=C(\overline I;C^2(\Gamma))\cap C^1(\overline I;C(\Gamma))
\end{equation}
with $h(t)$ admissible for every $t\in\overline I$ and let $V$ denote the normal velocity of $\Gamma[h(t)]$ in direction $n[h(t)]$.
Since $\partial_tX_h=\partial_th\,n_\Gamma$, its reference pullback is
\begin{equation}\label{eq:velocity_height}
 \widehat V_h:=V\circ X_h=\beta(h)\partial_th.
\end{equation}
Moreover, $\Theta_h=\mathrm{Id}+\mathcal Eh$ is affine in $h$, where
$(\mathcal Ev)(\Lambda(\gamma,\zeta))=\chi(\zeta)v(\gamma)n_\Gamma(\gamma)$ in $U_\Gamma$ and
$\mathcal Ev=0$ otherwise.
Hence $\partial_t\Theta_h=\mathcal E(\partial_th)$, and the transformed equations contain the transport coefficient
\begin{equation}\label{eq:geometric_transport}
 b_h=D\Theta_h^{-1}\mathcal E(\partial_th).
\end{equation}

Finally, for a family $\bh_\e=(h_{\e,k})_{k\in \indexI}$, we define
\[
 \Gamma_{\e,k}(t)=\e(k+\Gamma[h_{\e,k}(t)]),\qquad
 \Theta_\e(t,\e(k+y))=\e(k+\Theta_{h_{\e,k}(t)}(y)).
\]
Since $\Theta_h$ is the identity outside $U_\Gamma\Subset Y$, the map $\Theta_\e$ is the identity near every cell face, and we extend it by the identity outside the inclusion cells.
In particular, $\Theta_\e=\mathrm{Id}$ near $\partial\Omega$.
We identify the family $\bh_\e=(h_{\e,k})_{k\in\indexI}$ with the height field $h_\e\colon\overline{I}\times\Gamma_\e\to\R$ given by
\[
 h_\e(t,\e(k+\gamma))=h_{\e,k}(t,\gamma),\qquad \gamma\in\Gamma,\ k\in\indexI.
\]
The actual height on $\Gamma_\e$ is $\e h_\e$, whereas $h_\e$ is dimensionless and geometric quantities of $h_\e$ are understood cellwise, that is,
\[
 \beta(h_\e)(t,\gamma_{\e,k})=\beta(h_{\e,k}(t))(\gamma),\qquad
 \widehat H[h_\e](t,\gamma_{\e,k})=\widehat H[h_{\e,k}(t)](\gamma).
\]
with $\gamma_{\e,k}=\e(k+\gamma)$.
On $\Gamma_\e$, we write $\widehat n_\e=n_\e\circ\Theta_\e$, $\widehat H_\e=H_\e\circ\Theta_\e$, and $\widehat v_\e=v_\e\circ\Theta_\e$.
For $\gamma\in\Gamma$, the scaling identities are
\begin{equation*}
 \widehat n_\e(t,\gamma_{\e,k})=\widehat n_{h_{\e,k}(t)}(\gamma),\quad
 \widehat H_\e(t,\gamma_{\e,k})=\e^{-1}\widehat H[h_{\e,k}(t)](\gamma),\quad
 \widehat v_\e(t,\gamma_{\e,k})=\e\widehat V_{h_{\e,k}}(t,\gamma).
\end{equation*}
Finally, let $\widehat u_\e=u_\e\circ\Theta_\e$, $\widehat f_\e=f_\e\circ\Theta_\e$, $J_\e=\det D\Theta_\e$, $A_\e=J_\e D\Theta_\e^{-1}D\Theta_\e^{-T}$, and $b_\e=D\Theta_\e^{-1}\partial_t\Theta_\e$.
On the fixed physical reference domains, the transformed problem reads
\begin{empheq}[left={\systemlabel{sys:osmosis-fixed}{\widehat P_\e}\quad\empheqlbrace}]{alignat=2}
 \partial_t\widehat u_\e^{\pm}-b_\e\cdot\nabla\widehat u_\e^{\pm}
 -\frac{d_\e^{\pm}}{J_\e}\dive(A_\e\nabla\widehat u_\e^{\pm})
 &=\widehat f_\e^{\pm}
 &&\text{in }\Omega_\e^{\pm},\label{eq:fixed_bulk}\\
 d_\e^{\pm}N_\e\widehat u_\e^{\pm}
 +\e\widehat u_\e^{\pm}\beta(h_\e)\partial_th_\e&=0
 &&\text{on }\Gamma_\e,\label{eq:fixed_physical_boundary}\\
 \beta(h_\e)\partial_th_\e+\sigma \widehat H[h_\e]
 &=\lambda(\widehat u_{\e}^{-}-\widehat u_{\e}^{+})
 &&\text{on }\Gamma_\e,\label{eq:fixed_velocity}\\
 d^+\nabla\widehat u_\e^{+}\cdot\nu&=0 &&\text{on }\partial\Omega,\\
 \widehat u_\e^{\pm}(0)&=u_{\e,0}^{\pm},\qquad h_\e(0)=0.
 &&
\end{empheq}

\begin{remark}\label{rem:membrane_law}
\begin{itemize}
 \item Note that the curvature $\widehat H[h_\e]$ in \eqref{eq:fixed_velocity} is the mean curvature of the dimensionless cell interface $\Gamma[h_{\e,k}]\subset Y$, pulled back to $\Gamma$ and evaluated cellwise.
 It is not the curvature of the physical membrane $\Gamma_\e(t)$, which is $\widehat H_\e=\e^{-1}\widehat H[h_\e]$.
 Written entirely in terms of $\Gamma_\e$, the membrane law reads
 \[
 \beta(h_\e)\partial_t h_\e+\e\sigma\widehat H_\e=\lambda(\widehat u_\e^{-}-\widehat u_\e^{+})\qquad\text{on }\Gamma_\e.
 \]
 Since tangential derivatives on $\Gamma$ and on $\Gamma_\e$ are related by $\nabla_\Gamma=\e\nabla_{\Gamma_\e}$, the principal part of $\e\sigma\widehat H_\e$ is $-\e^2\sigma\Delta_{\Gamma_\e}h_\e$.
 On the physical membrane, the height equation is therefore a nonlinear parabolic equation whose diffusion coefficient degenerates like $\e^2$.
 On the fixed reference membrane $\Gamma$, the same equation is uniformly parabolic with coefficients independent of $\e$.
 \item The time-dependent Jacobian also gives the conservative formulation
\begin{equation}\label{eq:fixed_conservative}
 \partial_t(J_\e\widehat u_\e^{\pm})
 -\dive\bigl(d_\e^{\pm} A_\e\nabla\widehat u_\e^{\pm}
                    +J_\e b_\e\widehat u_\e^{\pm}\bigr)
 =J_\e\widehat f_\e^{\pm}
\end{equation}
via $\partial_tJ_\e=\dive(J_\e b_\e)$ (by Jacobi's formula and the Piola identity).
\end{itemize}
\end{remark}
For $\tau\in(0,T^*)$, $I_\tau=(0,\tau)$, and $\beta\in(0,1)$, we introduce the Hölder-type space
\begin{align*}
  \mathcal X^\beta(I_\tau)&:=\mathcal X(I_\tau)\cap C(\overline{I_\tau};C^{2+\beta}(\Gamma)),
\end{align*}
with $\mathcal X(I_\tau)$ defined in \eqref{eq:defintion_XI}.
In addition, we write $\mathcal X_\e(I_\tau)$ and $\mathcal X_\e^\beta(I_\tau)$ for the corresponding products over $k\in\indexI$, with the maximum norm.
For either phase, let
\[
 Q_\e^\pm(I_\tau):=\{(t,x):t\in I_\tau,\ x\in\Omega_\e^\pm(t)\}.
\]

\begin{definition}[Weak microscopic solution]\label{def:weak_solution}
Let $\tau\in(0, T^*)$.
The functions $(u_\e^\pm,\bh_\e)$ are a weak solution of \systemref{sys:osmosis}{P_\e} on $I_\tau$ if the following holds:
\begin{enumerate}
 \item The heights satisfy $\bh_\e\in\mathcal X_\e(I_\tau)$ with $\bh_\e(0)=0$ and
 \[
 \max_{k\in\indexI}\sup_{t\in\overline{I_\tau}}\|h_{\e,k}(t)\|_{C^1(\Gamma)}\le\eta^*.
 \]
 In particular, each $\Theta_\e(t)$ is a diffeomorphism with $\Omega_\e^\pm(t)=\Theta_\e(t,\Omega_\e^\pm)$ and $\Gamma_\e(t)=\Theta_\e(t,\Gamma_\e)$.
 \item The concentrations $u_\e^\pm$ are defined on $Q_\e^\pm(I_\tau)$, and their pullbacks satisfy
 \[
 \widehat u_\e^\pm=u_\e^\pm\circ\Theta_\e\in C(\overline{I_\tau};L^2(\Omega_\e^\pm))\cap L^2(I_\tau;H^1(\Omega_\e^\pm))\cap L^\infty(I_\tau\times\Omega_\e^\pm).
 \]
 \item The bulk equations and the impermeability conditions hold in the weak conservative sense: for every $\varphi\in C^\infty(\overline{I_\tau}\times\overline{\Omega_\e^\pm})$ with $\varphi(\tau)=0$,
 \begin{multline}\label{eq:weak_form}
  -\int_0^\tau\!\int_{\Omega_\e^\pm}J_\e\widehat u_\e^\pm\,\partial_t\varphi\di x\di t
  +\int_0^\tau\!\int_{\Omega_\e^\pm}\bigl(d_\e^\pm A_\e\nabla\widehat u_\e^\pm+J_\e b_\e\widehat u_\e^\pm\bigr)\cdot\nabla\varphi\di x\di t\\
  =\int_0^\tau\!\int_{\Omega_\e^\pm}J_\e\widehat f_\e^\pm\varphi\di x\di t
  +\int_{\Omega_\e^\pm}u_{\e,0}^\pm\varphi(0)\di x.
 \end{multline}
 \item The membrane law \eqref{eq:fixed_velocity} holds almost everywhere on $I_\tau\times\Gamma_\e$, with the traces of $\widehat u_\e^\pm$ on $\Gamma_\e$.
\end{enumerate}
\end{definition}

\subsection{Main results}
\label{ssec:main_results}
We state the principal conclusions of this work.
Uniform existence requires only Assumptions~\ref{ass:geometry}--\ref{ass:sources}, while homogenization additionally relies on Assumptions~\ref{ass:macro_initial} and~\ref{ass:macro_sources}.

\paragraph{A common existence interval}
The first main result provides weak solutions on a time interval independent of $\e$; it is proved in Section~\ref{sec:analysis}.

\begin{theorem}[Uniform local existence]
\label{thm:prepared_schauder}
There exist $\e$-independent parameters $0<\beta<\alpha<1$, $M,B>0$, and
$T\in(0,T^*]$ such that Problem~$\systemref{sys:osmosis}{P_\e}$ admits a weak
solution $(u_\e^\pm,\bh_\e)$ in the sense of \cref{def:weak_solution} with
$\bh_\e\in\mathcal X_\e^\beta(I)$ on $I:=(0,T)$.
Moreover,
\begin{equation}\label{eq:prepared_uniform_result}
    \max_{k\in\indexI}
    \|h_{\e,k}\|_{C(\overline{I};C^{2+\beta}(\Gamma))}
    \le \frac{\eta^*}2,
    \qquad
    \max_{k\in\indexI}
    \|\partial_t h_{\e,k}\|_{C(\overline{I}\times\Gamma)}
    \le M,
\end{equation}
and
\[
    \max_{\pm}
    \|u_\e^\pm\|_
    {L^\infty(Q_\e^\pm)}
    \le u^*+1,\qquad 
    \max_{\pm}
    \sup_{t\in\overline{I}}
    \bigl\|
        \widehat u_{\e}^\pm(t)
    \bigr\|_{C^\alpha(\Gamma_\e)}
    \le \e^{-\alpha}B.
\]
Here, $u^*$ is the bound from Assumption~\ref{ass:initial_data}.
\end{theorem}

\paragraph{Homogenization}
For an admissible height field $h=h(t,x)$ with $\|h(t,x)\|_{C^1(\Gamma)}\le\eta^*$, we define the phase volumes
\begin{equation}\label{eq:phase_volumes}
 \theta^{\pm}(t,x)=\int_{Y^{\pm}}J_{h(t,x)}(y)\di y
              =|Y^{\pm}[h(t,x)]|.
\end{equation}
Let $H^1_{\#}(Y^{+})$ denote the space with matching traces on opposite cell faces and put
\[
 \mathcal V=\left\{\psi\in H^1_{\#}(Y^{+}):
                              \int_{Y^{+}}\psi\di y=0\right\}.
\]
For each admissible height $h$ and $j=1,2,3$, let $\chi_j[h]\in\mathcal V$ be the unique solution of
\begin{equation}\label{eq:corrector_basis}
 \int_{Y^{+}}A_h(e_j+\nabla_y\chi_j[h])\cdot\nabla_y\psi\di y=0
              \qquad(\psi\in\mathcal V).
\end{equation}
We introduce the effective macroscopic diffusivity
\begin{equation}\label{eq:effective_tensor}
 D_{ij}(h)=d^{+}\int_{Y^{+}}e_i\cdot A_h(e_j+\nabla_y\chi_j[h])\di y.
\end{equation}
The existence of these correctors and the properties of $D(h)$ are proved in \cref{lem:cell_correctors}.
We introduce the limit problem
\begin{empheq}[left={\systemlabel{sys:homogenized-fixed}{\widehat P_0}\quad\empheqlbrace}]{alignat*=2}
 \partial_t(\theta^{+}(h)u^{+})-\dive_x(D(h)\nabla_xu^{+})
 &=\theta^{+}(h)f^{+}
 &&\quad\text{in }I\times\Omega,\\
 \partial_t(J_h\widehat u^{-})
 -\dive_y\bigl(d^{-}A_h\nabla_y\widehat u^{-}+J_hb_h\widehat u^{-}\bigr)
 &=J_h\widehat f^{-}
 &&\quad\text{in }I\times\Omega\times Y^{-},\\
 \bigl(d^{-}A_h\nabla_y\widehat u^{-}+J_hb_h\widehat u^{-}\bigr)\cdot n_\Gamma
 &=0
 &&\quad\text{on }I\times\Omega\times\Gamma,\\
 \beta(h)\partial_th+\sigma\widehat H[h]
 &=\lambda(\widehat u^--u^+)
 &&\quad\text{on }I\times\Omega\times\Gamma,\\
 D(h)\nabla_xu^{+}\cdot\nu
 &=0
 &&\quad\text{on }I\times\partial\Omega,\\
 u^{+}(0)=u_0^+,\quad \widehat u^{-}(0)&=u_0^-,\quad h(0)=0,
 &&
\end{empheq}
where the transformation quantities $J_h,A_h,b_h,\beta(h)$ are defined as before and $\widehat f^-(t,x,y)=f^-(t,x,\Theta_{h(t,x)}(y))$.

\begin{definition}[Solution of the homogenized problem]
\label{def:weak_homogenized}
A triple $(u^+,\widehat u^-,h)$ is a weak solution of \systemref{sys:homogenized-fixed}{\widehat P_0} on $I$ if the following conditions hold.
\begin{enumerate}[label=(\alph*)]
 \item \emph{Regularity and admissibility.}
 \begin{align*}
  u^+&\in L^2(I;H^1(\Omega))\cap L^\infty(I\times\Omega),\\
  \widehat u^-&\in C(\overline I;L^2(\Omega\times Y^-))\cap L^2\big(I\times \Omega;H^1(Y^-)\big)\cap L^\infty(I\times\Omega\times Y^-),\\
  h&\in L^\infty\big(\Omega;C(\overline{I};C^{2}(\Gamma))\cap W^{1,\infty}(I;C(\Gamma))\big),
 \end{align*}
 and $h$ is admissible, i.e.\ $\|h(t,x)\|_{C^1(\Gamma)}\le\eta^*$, for all $t\in\overline I$ and a.a.~$x\in\Omega$.

 \item \emph{Macroscopic equation.}
 For every $\varphi\in C^\infty(\overline{I}\times\overline\Omega)$ with $\varphi(T)=0$,
 \begin{multline*}
  -\int_0^T\int_\Omega\theta^+(h) u^+\partial_t\varphi\di x\di t
  +\int_0^T\int_\Omega D(h)\nabla_xu^+\cdot\nabla_x\varphi\di x\di t\\
  =\int_0^T\int_\Omega\theta^+(h)f^+\varphi\di x\di t
  +|Y^+|\int_\Omega u_0^+(x)\varphi(0,x)\di x.
 \end{multline*}

 \item \emph{Cell problems.}
 \begin{itemize}
  \item[$(i)$] For almost every $x\in\Omega$ and for every $\phi\in C^\infty(\overline{I}\times\overline{Y^-})$ with $\phi(T)=0$,
  \begin{multline*}
   -\int_0^T\int_{Y^-}J_h\widehat u^-\,\partial_t\phi\di y\di t
   +\int_0^T\int_{Y^-}\bigl(d^-A_h\nabla_y\widehat u^-+J_h b_h \widehat u^-\bigr)\cdot\nabla_y\phi\di y\di t\\
   =\int_0^T\int_{Y^-}J_h\,(f^-\circ\Theta_h)\,\phi\di y\di t
   +\int_{Y^-}u_0^-\phi(0)\di y.
  \end{multline*}
  \item[$(ii)$] For almost every $x\in\Omega$, $h(0,x,\cdot)=0$ and
  \[
   \beta(h)\partial_th+\sigma\widehat H[h]=\lambda\bigl(\widehat u^--u^+\bigr)
   \qquad\text{on }\Gamma,\ \text{for a.e. } t\in I.
  \]
 \end{itemize}
\end{enumerate}
\end{definition}

\begin{theorem}[Homogenization]
\label{thm:homogenization}
Assume \ref{ass:geometry}--\ref{ass:sources}, \ref{ass:macro_initial}, and \ref{ass:macro_sources}, and let $(u_\e^\pm,\bh_\e)$ be weak solutions as in \cref{thm:prepared_schauder}.
Then every sequence $\e\downarrow0$ has a subsequence, not relabeled, and limits
\[
 u^+\in L^2(I;H^1(\Omega))\cap L^\infty(I\times\Omega),\qquad
 \widehat u^-\in L^\infty\big(\Omega;C(\overline{I}\times\overline{Y^-})\cap L^2(I;H^1(Y^-))\big),\]
 and
 \[
 h\in L^\infty\big(\Omega;C(\overline{I};C^{2+\beta'}(\Gamma))\cap W^{1,\infty}(I;C(\Gamma))\big),
\]
such that, for every $1\le q<\infty$ and $0<\beta'<\beta$, the following holds:
\begin{enumerate}
 \item[(i)] \underline{Convergence:}
 It holds $P_\e^+\widehat u_\e^+\to u^+$ strongly in $L^q(I\times\Omega)$ and weakly in $L^2(I;H^1(\Omega))$.
 Moreover,
 \[
  \mathcal T_\e\widehat u_\e^-\to \widehat u^- \ \text{in } L^q(\Omega;C(\overline{I}\times\overline{Y^-})),\quad
  \mathcal R_\e(h_{\e,k})_k\to h \ \text{in } L^q(\Omega;C(\overline{I};C^{2+\beta'}(\Gamma))).
 \]
 The unfolding operators $\mathcal T_\e$ and $\mathcal R_\e$ are given in Appendix~\ref{app:unfolding}.
 \item[(ii)] \underline{Limit problem:}  $(u^+,\widehat u^-,h)$ is a weak solution of \systemref{sys:homogenized-fixed}{\widehat P_0} in the sense of \cref{def:weak_homogenized}.
 \item[(iii)] \underline{Bounds:} The limit inherits the uniform bounds of \cref{thm:prepared_schauder}: $\|u^+\|_{L^\infty}\le u^*+1$, $\|\widehat u^-\|_{L^\infty}\le u^*+1$, and, for almost every $x\in\Omega$,
 \[
  \sup_{t\in\overline{I}}\|h(t,x)\|_{C^{2+\beta}(\Gamma)}\le \frac{\eta^*}2,
  \qquad
  \|\partial_th(\cdot,x)\|_{L^\infty(I\times\Gamma)}\le M .
 \]
\end{enumerate}
\end{theorem}

\paragraph{The limit problem in moving geometry} For almost every $x\in\Omega$, define
\begin{equation}\label{eq:limit_moving_domains}
 Y^{\pm}(t,x)=Y^{\pm}[h(t,x)],\qquad \Gamma(t,x)=\Gamma[h(t,x)],
 \qquad \theta^{\pm}(t,x)=|Y^{\pm}(t,x)|,
\end{equation}
so that $\theta^++\theta^-=1$, and the physical inner concentration
\begin{equation}\label{eq:limit_physical_concentration}
 u^{-}(t,x,\cdot)
 =\widehat u^{-}(t,x,\cdot)\circ\Theta_{h(t,x)}^{-1}
 \qquad\text{on } Y^{-}(t,x).
\end{equation}
By \cref{lem:app_geometry_flux}, applied for almost every $x$, \systemref{sys:homogenized-fixed}{\widehat P_0} is the conservative weak form of the moving-cell system \systemref{sys:homogenized-moving}{P_0} stated in the introduction, completed by $D\nabla_xu^+\cdot\nu=0$ on $\partial\Omega$ and the initial conditions $u^+(0)=u_0^+$, $u^-(0)=u_0^-$, $\Gamma(0,x)=\Gamma$.
There, $V$ is the normal velocity with $V\circ X_h=\beta(h)\partial_th$.
By \cref{rem:physical_tensor}, $D(h)$ coincides with the classical Neumann-homogenized tensor of the physical outer cell $Y^+(t,x)$.

\begin{remark}[Interpretation of the limit problem]\label{rem:homogenization}
\begin{enumerate}[label=(\arabic*)]
 \item \emph{Coupling through geometry.}
As in the $\e$-problem, there is no solute exchange between the phases.
The phases interact only through the membrane motion, which enters the outer
equation via the storage factor $\theta^+(h)$ and the effective diffusivity $D(h)$.
While $\theta^+(h)$ only sees the volume of the inclusion, $D(h)$ also responds to
the shape of $Y^-(t,x)$.\footnote{Just as a simple example, rotating an elliptic inclusion leaves $\theta^+(h)$ unchanged but rotates the principal axes of $D(h)$.}
Such a feedback through effective porosity and diffusivity also occurs for evolving
solid microstructures \cite{GahnPop2023,WiedemannPeter2023,GahnPeterPopWiedemann2024},
where, however, the geometry is described by finitely many parameters such as a radius
or a layer thickness.
 \item \emph{Consistency.}
 If the membrane motion is frozen, the limit reduces to the classical perforated-domain equation with storage $\theta^+(h)$; e.g., for a ball $Y^-$, $f^\pm=0$, $u_0^+\equiv\overline u$, and $u_0^-\equiv\overline u+\frac{\sigma}{\lambda}H_\Gamma$, one has $h\equiv0$ and $D(0)$ is the classical tensor for impermeable holes.
 \item \emph{Uniqueness.}
 We expect weak--strong uniqueness: if $\nabla_xu_2^+\in L^\infty(I\times\Omega)$ and
 $\partial_tu_2^+\in L^2(I;L^3(\Omega))$, \cref{thm:cell_stability} controls $h_1-h_2$ by  $u_1^+-u_2^+$ for a.e.\ $x$, and the Lipschitz dependence of $D(h)$ and $\theta^+(h)$ on $h$ suggests a Gronwall argument.
\end{enumerate}
\end{remark}

\section{Analysis of the \texorpdfstring{$\e$}{epsilon}-problem}\label{sec:analysis}

Throughout this section we impose Assumptions~\ref{ass:geometry}--\ref{ass:sources} and prove \cref{thm:prepared_schauder}.
On the scale of a single inclusion, the outer diffusion is of order $\e^{-2}$, so the outer concentration has no uniform modulus of continuity in time; the membrane law, however, only sees the concentration traces on each interface.
We localize the solution around each inclusion (\cref{ssec:rescaling}), derive spatial trace bounds that are uniform in $\e$ for prescribed heights (\cref{ssec:prescribed_motion}), show that spatial control alone keeps the membrane evolution bounded on a short interval (\cref{ssec:schauder}), and combine both by Schauder's fixed point theorem (\cref{ssec:proof_wellposed}).

\subsection{Localization and cellwise rescaling}
\label{ssec:rescaling}
\begin{figure}[htbp]
    \centering
    \includegraphics[width=.4\linewidth]{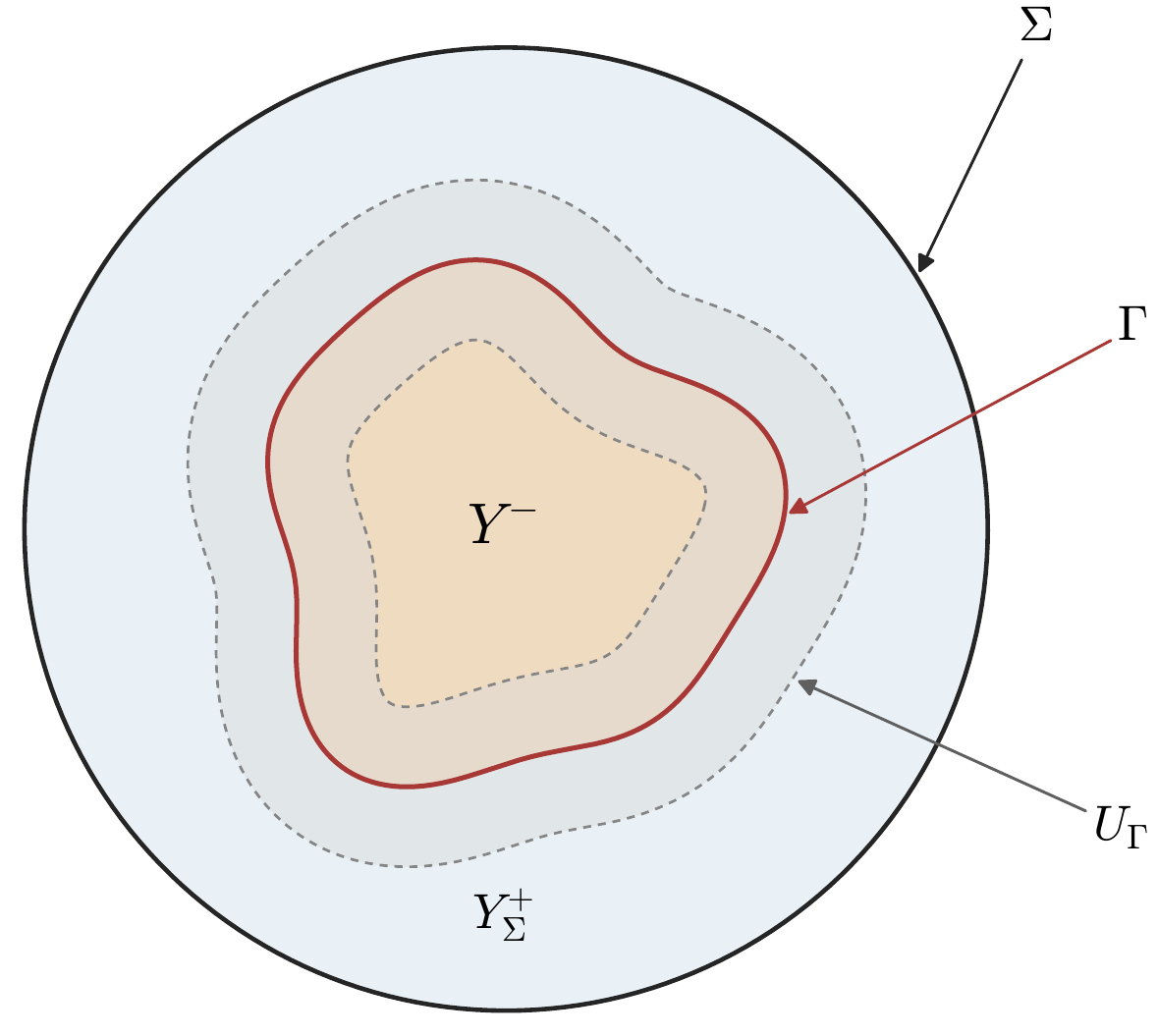}
    \caption{Reference geometry $Y_\Sigma$ of the localized cell problem: inner phase $Y^-$, local outer collar $Y_\Sigma^+$, membrane $\Gamma$, and tubular neighborhood $U_\Gamma$ in which the Hanzawa transformation acts.}
    \label{fig:cellwise_rescaling}
\end{figure}

The membrane law is driven by the concentration traces on the individual interfaces $\Gamma_{\e,k}(t)$.
To estimate these traces uniformly in $\e$ and $k$, we rescale a fixed neighborhood of each inclusion to the reference cell.
To that end, we choose a smooth domain $Y_\Sigma\Subset Y$ with $\overline{Y^-}\cup\overline{U_\Gamma}\subset Y_\Sigma$ and set
\[
 \Sigma=\partial Y_\Sigma,\qquad
 Y_\Sigma^{+}=Y_\Sigma\setminus\overline{Y^{-}},
\]
see \cref{fig:cellwise_rescaling}.
Since $\Theta_h$ is the identity outside $U_\Gamma$, it maps $Y^-$ onto $Y^-[h]$ and $Y_\Sigma^+$ onto $Y_\Sigma^+[h]=Y_\Sigma\setminus\overline{Y^-[h]}$, and it equals the identity near $\Sigma$ for every admissible height.
For $k\in\indexI$ and $t\in I_\tau$ we define the cellwise restrictions
\[
 \widehat\mu_{\e,k}^{\pm}(t,y)=\widehat u_\e^{\pm}(t,\e(k+y)),\qquad
 \widehat f_{\e,k}^{\pm}(t,y)=\widehat f_\e^{\pm}(t,\e(k+y)),
\]
for $y\in Y^-$ and $y\in Y_\Sigma^+$, respectively, and write $\Theta_{\e,k}(t)=\Theta_{h_{\e,k}(t)}$.
Using the scalings $\nabla_x=\e^{-1}\nabla_y$, $b_\e=\e b_h$, $N_\e=\e^{-1}N_h$, and $\widehat v_\e=\e\beta(h_{\e,k})\partial_t h_{\e,k}$ from \cref{ssec:fixed_domain}, the restrictions satisfy, in the weak sense,
\begin{empheq}[
    left={\empheqlbrace}
]{alignat=2}
    \e^2\partial_t\widehat\mu_{\e,k}^{+}
    -\e^2 b_{h_{\e,k}}\cdot\nabla_y\widehat\mu_{\e,k}^{+}
    -d^{+}J_{h_{\e,k}}^{-1}
    \dive_y(A_{h_{\e,k}}\nabla_y\widehat\mu_{\e,k}^{+})
    &=\e^2\widehat f_{\e,k}^{+}
    &&\quad\text{in }Y_\Sigma^+,
    \label{eq:fixed_cell:plus}\\
    \partial_t\widehat\mu_{\e,k}^{-}
    -b_{h_{\e,k}}\cdot\nabla_y\widehat\mu_{\e,k}^{-}
    -d^{-}J_{h_{\e,k}}^{-1}
    \dive_y(A_{h_{\e,k}}\nabla_y\widehat\mu_{\e,k}^{-})
    &=\widehat f_{\e,k}^{-}
    &&\quad\text{in }Y^-,
    \label{eq:fixed_cell:minus}\\
    d^{+}N_{h_{\e,k}}\widehat\mu_{\e,k}^{+}
    +\e^2\widehat\mu_{\e,k}^{+}
    \beta(h_{\e,k})\partial_t h_{\e,k}
    &=0
    &&\quad\text{on }\Gamma,
    \label{eq:fixed_cell:flux:plus}\\
    d^{-}N_{h_{\e,k}}\widehat\mu_{\e,k}^{-}
    +\widehat\mu_{\e,k}^{-}
    \beta(h_{\e,k})\partial_t h_{\e,k}
    &=0
    &&\quad\text{on }\Gamma,
    \label{eq:fixed_cell:flux:minus}\\
    \beta(h_{\e,k})\partial_t h_{\e,k}
    +\sigma \widehat H[h_{\e,k}]
    &=\lambda\jump{\widehat\mu_{\e,k}}
    &&\quad\text{on }\Gamma,
    \label{eq:fixed_cell:velocity}
\end{empheq}
Here \eqref{eq:fixed_cell:plus} and \eqref{eq:fixed_cell:flux:plus} hold for test functions vanishing near $\Sigma$, and the initial values are the rescaled initial data.
No condition on $\Sigma$ is imposed or needed.
The estimates below are local near $\Gamma$ and use only the global $L^\infty$ bound on the collar, which has positive distance from $\Sigma$.
 
The factor $\e^2$ in front of the time derivative in \eqref{eq:fixed_cell:plus} is the only difference between the phases; it is treated by the fast time $s=t/\e^2$ in \cref{pro:moving_holder}.

\subsection{Bulk diffusion with prescribed interface motion}
\label{ssec:prescribed_motion}
In this subsection, we derive a priori estimates for solutions of Problem \systemref{sys:osmosis-fixed}{\widehat P_\e} that are uniform with respect to $\e$ for a class of prescribed height functions.
For a family $\bh=(h_k)\in\mathcal X_\e(I_\tau)$, we introduce the geometric bounds
\begin{equation}\tag{\ensuremath{\mathrm{GB}_{M}}}\label{ass:geometry_bootstrap}
 \max_k\|h_{k}\|_{C(\overline{I_\tau};C^2(\Gamma))}\le\frac{\eta^*}2,\qquad
 \max_k\|\partial_th_{k}\|_{L^\infty(I_\tau\times\Gamma)}\le M
\end{equation}
and the closed subset
\[
    \mathcal H_{M}(I_\tau)
    :=
    \bigl\{
        \bh\in\mathcal X_\e(I_\tau):
        \bh\text{ satisfies }\eqref{ass:geometry_bootstrap}
    \bigr\}
\]
the corresponding class of admissible interface motions.
Here $M>0$ is a constant later to be fixed and $\eta^*$ is the geometric
threshold fixed in \cref{ssec:fixed_domain} (see also \cref{lem:geometric_estimates}).
Since
\[
\max_k\|h_{k}\|_{C(\overline{I_\tau};C^1(\Gamma))}\le\max_k\|h_{k}\|_{C(\overline{I_\tau};C^2(\Gamma))}\le\frac{\eta^*}2,
\]
the corresponding Hanzawa transforms are well-defined diffeomorphisms.
The factor $1/2$ leaves room for the smooth approximations in \cref{lem:app_motion_approximation}.\footnote{Hence, any bound $q\eta^*$ with $q<1$ would be sufficient here.}
In the following, we consider Problem \systemref{sys:osmosis-fixed}{\widehat P_\e} for prescribed $\bh\in\mathcal H_{M}(I_\tau)$ with $\bh(0)=0$.
On $\{\|h\|_{C^2(\Gamma)}\le\eta^*/2\}$, we set
\[
  b_*:=\inf\beta(h)>0,\qquad H^*:=\sup\|\widehat H[h]\|_{L^\infty(\Gamma)}<\infty,
\]
which depend only on the reference geometry and are now fixed.

\begin{lemma}
\label{lem:prescribed_bulk}
Let $I_\tau=(0,\tau)\subset I^*$ and let $\bh\in\mathcal H_{M}(I_\tau)$ with $\bh(0)=0$.
Problem~\systemref{sys:osmosis-fixed}{\widehat P_\e}, with the membrane motion prescribed by $\bh$, admits a unique weak solution
\[
    \widehat u_\e^\pm
    \in C(\overline{I_\tau};L^2(\Omega_\e^\pm))
       \cap L^2(I_\tau;H^1(\Omega_\e^\pm)).
\]
In addition, the reference mass densities \(q_\e^\pm:=J_\e\widehat u_\e^\pm\) satisfy
\[
    q_\e^\pm
    \in L^2(I_\tau;H^1(\Omega_\e^\pm))
       \cap H^1(I_\tau;H^1(\Omega_\e^\pm)').
\]
\end{lemma}

\begin{proof}
For a prescribed motion the two phases decouple.
Fix one phase, write $D=\Omega_\e^\pm$ and $d=d_\e^\pm$, and set $K=dA_\e/J_\e$ and $c=b_\e-K\nabla\log J_\e$.
Then $q=J_\e\widehat u_\e^\pm$ solves the conservative problem
\begin{equation}\label{eq:prescribed_density}
\begin{aligned}
    \partial_t q-\dive(K\nabla q+cq)&=J_\e\widehat f
    &&\text{in }I_\tau\times D,\\
    (K\nabla q+cq)\cdot\nu_D&=0
    &&\text{on }I_\tau\times\partial D,\qquad q(0)=u_{\e,0}^\pm,
\end{aligned}
\end{equation}
where the conormal condition encodes impermeability on $\Gamma_\e$ and, for $D=\Omega_\e^+$, the no-flux condition on $\partial\Omega$, where $c=0$.
For fixed $\e$, \cref{lem:geometric_estimates} shows that $K$ is bounded and uniformly elliptic, that $c$ is bounded, and that $J_\e^{\pm1}\in C(\overline{I_\tau};W^{1,\infty}(D))$.
Standard variational theory \cite[Chapter~XVIII, Sections~3 and 4.2]{DautrayLions1992} gives a unique $q\in L^2(I_\tau;H^1(D))\cap H^1(I_\tau;H^1(D)')\hookrightarrow C(\overline{I_\tau};L^2(D))$, and $\widehat u_\e^\pm=J_\e^{-1}q$ has the stated regularity.
\end{proof}

By the Hanzawa transformation, the solutions $\widehat u_\e^\pm$ from \cref{lem:prescribed_bulk} define functions on the moving domains via
\(
    u_\e^\pm(t)
    :=\widehat u_\e^\pm(t)\circ\Theta_\e(t)^{-1}.
\)
We now derive energy estimates and nonnegativity for these solutions in the moving domains.

\begin{lemma}
\label{pro:basic_estimates}
Let $I_\tau=(0,\tau)\subset I^*$ and $\bh_\e\in\mathcal H_{M}(I_\tau)$ satisfy $\bh_\e(0)=0$.
Let $u_\e^\pm$ be the corresponding solutions via \cref{lem:prescribed_bulk}.
\begin{itemize}
\item[$(i)$] The concentrations satisfy the energy estimate
    \begin{multline*}
        \sup_{t\in I_\tau}
        \sum_{\pm}
        \|u_\e^{\pm}(t)\|_{L^2(\Omega_\e^{\pm}(t))}^2
        +
        \int_0^\tau
        \Bigl(
            d^{+}\|\nabla u_\e^{+}\|_{L^2(\Omega_\e^{+}(t))}^2
            +
            \e^2 d^{-}
            \|\nabla u_\e^{-}\|_{L^2(\Omega_\e^{-}(t))}^2
        \Bigr)\,\mathrm dt
        \\
        \lesssim_{M}\ 
        \sum_{\pm}
        \left(
            \|u_{\e,0}^{\pm}\|_{L^2(\Omega_\e^{\pm})}^2
            +
            \|f_\e^{\pm}\|_{L^2(I_\tau\times\Omega)}^2
        \right).
    \end{multline*}

\item[$(ii)$] \emph{Comparison.} Let $q^\pm$ be functions on the moving phases with the regularity of $u_\e^\pm$ that satisfy, in the weak sense,
\begin{gather*}
 \partial_tq^\pm-d_\e^\pm\Delta q^\pm\ge0\ \text{ in }\Omega_\e^\pm(t),\qquad
 \nabla q^+\cdot\nu=0\ \text{ on }\partial\Omega,\\
 \mp\bigl(d_\e^\pm\nabla q^\pm\cdot n_\e+v_\e q^\pm\bigr)\ge0\ \text{ on }\Gamma_\e(t),
\end{gather*}
and $q^\pm(0)\ge0$. Then $q^\pm\ge0$.
In particular, $u_\e^\pm\ge0$ whenever $u_{\e,0}^{\pm}\ge0$ and $f_\e^{\pm}\ge0$.
\end{itemize}
\end{lemma}
\begin{proof}
$(i)$.
For sufficiently regular solutions, testing the bulk equations with $u_\e^\pm$ and using Reynolds' formula yields
\begin{equation*}
 \frac12\frac{\mathrm d}{\mathrm dt}
 \sum_{\pm}\int_{\Omega_\e^{\pm}(t)}|u_\e^{\pm}|^2
 +\sum_{\pm}d_\e^{\pm}
 \int_{\Omega_\e^{\pm}(t)}|\nabla u_\e^{\pm}|^2
 =
 \sum_{\pm}\int_{\Omega_\e^{\pm}(t)}f_\e^{\pm}u_\e^{\pm}
 +\frac12\int_{\Gamma_\e(t)}
 \bigl(|u_\e^+|^2-|u_\e^-|^2\bigr)v_\e.
\end{equation*}
We apply \cref{lem:scaled_trace} to the last term and choose $a$ so that
$CMa\le\frac12\min\{d^{+},d^{-}\}$, using $|v_\e|\le\e CM$.
Since $\e\le1$, the resulting gradient terms can be absorbed into both diffusion contributions.
The claim then follows via Young's inequality and Gronwall's lemma.

$(ii)$.
Testing with $-\operatorname{neg}(q^\pm)$ and using Reynolds' formula gives
\begin{equation*}
 \frac12\frac{\mathrm d}{\mathrm dt}
       \|\operatorname{neg}(q^{\pm})\|_{L^2(\Omega_\e^{\pm}(t))}^2
 +d_\e^{\pm}\|\nabla\operatorname{neg}(q^{\pm})\|_{L^2(\Omega_\e^{\pm}(t))}^2
 \le\frac12\int_{\Gamma_\e(t)}
       |v_\e||\operatorname{neg}(q^{\pm})|^2\di S.
\end{equation*}
The boundary term is absorbed exactly as in $(i)$.
Since the initial negative part vanishes, Gronwall's lemma gives $q^\pm\ge0$.
For $q^\pm=u_\e^\pm$, the source contribution is nonnegative when $f_\e^\pm\ge0$.
Both arguments extend to weak solutions by approximation (\cref{lem:app_motion_approximation,lem:app_bulk_continuity}).
\end{proof}

\begin{lemma}
\label{lem:short_sup}
Let $0<\tau\le\min\{T^*,T_{\mathrm{geom}}\}$ and $\bh\in\mathcal H_{M}(I_\tau)$ with $\bh(0)={\bf0}$ and let $u_\e^{\pm}$ be the corresponding solution.
Here $T_{\mathrm{geom}}>0$ depends only on the local geometry of $\Gamma$ inside $Y$. It holds
\begin{equation}\label{eq:short_sup}
 \max_{\pm}\|u_\e^{\pm}\|_{L^\infty(Q_\e^\pm)}
 \le (u^*+\tau F^*)
 \exp\!\left(C\bigl((1+M)\sqrt \tau+(1+M)^2\tau\bigr)\right).
\end{equation}
\end{lemma}
\begin{proof}
Since the membrane velocity gives the Robin term an indefinite sign, constant barriers do not suffice.
As for oblique derivative problems, we multiply by an exponential weight that reverses the sign of the Robin term.
The weight is localized in a collar of width $\ell$ around the membranes, and $\ell$ is chosen depending on $\tau$ so that the bound tends to $u^*$ as $\tau\downarrow0$, uniformly in $\e$.

\emph{Weights.}
Let $\ell_0\in(0,1]$ be such that $\operatorname{dist}(\cdot,\Gamma)$ is $C^3$ on both sides of $\Gamma$ in collars of width $\ell_0$ that do not meet $\Sigma$, and let $\chi_0\in C^\infty([0,\infty))$ be nonnegative with $\chi_0'(0)=-1$ and $\chi_0=0$ on $[1,\infty)$.
For $0<\ell\le\ell_0$ put $\psi_\ell^\pm=\ell\chi_0(\operatorname{dist}(\cdot,\Gamma)/\ell)$ in the collar of $Y^\pm$, extended by zero.
Then $0\le\psi_\ell^\pm\le C\ell$, $|\nabla\psi_\ell^\pm|\le C$, and $|D^2\psi_\ell^\pm|\le C/\ell$, and $\nabla\psi_\ell^\pm$ is the outward normal of the respective phase on $\Gamma$.
By \eqref{ass:geometry_bootstrap} and \cref{lem:geometric_estimates}, $N_h\psi_\ell^-\ge c$ and $-N_h\psi_\ell^+\ge c>0$ on $\Gamma$.
Let $W_\ell^-=\exp(\omega\psi_\ell^-)$ and $W_\ell^+=\exp(\e^2\omega\psi_\ell^+)$, and let $W_{\e,\ell}^\pm$ be their cellwise rescalings, with $W_{\e,\ell}^+=1$ outside the collars $\e(k+Y_\Sigma^+)$.
Then $1\le W_{\e,\ell}^\pm\le e^{C\omega\ell}$, and $\nabla W_{\e,\ell}^+\cdot\nu=0$ on $\partial\Omega$.

\emph{Robin and bulk inequalities.}
With $\widehat V_k=\beta(h_k)\partial_th_k$, so that $|\widehat V_k|\le CM$, the chain rule gives
\begin{align*}
 d^-N_hW_\ell^-+\widehat V_kW_\ell^-&\ge W_\ell^-(d^-\omega c-CM),\\
 -d^+N_hW_\ell^+-\e^2\widehat V_kW_\ell^+&\ge\e^2W_\ell^+(d^+\omega c-CM),
\end{align*}
and both right-hand sides are nonnegative for $\omega=C(1+M)$.
Since $N_\e=\e^{-1}N_h$ and $\widehat v_\e=\e\widehat V_k$, this yields $\mp\bigl(d_\e^\pm N_\e W_{\e,\ell}^\pm+\widehat v_\e W_{\e,\ell}^\pm\bigr)\ge0$ on $\Gamma_\e$.
For the bulk operator $\mathcal L_\e^\pm z=b_\e\cdot\nabla z+d_\e^\pm J_\e^{-1}\dive(A_\e\nabla z)$, direct differentiation in the rescaled variables, \cref{lem:geometric_estimates}, and the bounds on $\psi_\ell^\pm$ give
\[
 |\mathcal L_\e^\pm W_{\e,\ell}^\pm|\le K_\ell W_{\e,\ell}^\pm,\qquad
 K_\ell=C\bigl((1+M)/\ell+(1+M)^2\bigr),
\]
uniformly in $\e\le1$.
In the outer phase, the factor $\e^2$ in the exponent of $W_\ell^+$ exactly compensates the fast diffusion $\e^{-2}d^+$ on the cell scale.

\emph{Comparison.}
Let $A_\ell(t)=(u^*+tF^*)e^{K_\ell t}$ and $\widehat Z_\e^\pm=A_\ell W_{\e,\ell}^\pm$.
Since $A_\ell'-K_\ell A_\ell\ge F^*\ge|\widehat f_\e^\pm|$ by \ref{ass:sources}, $\widehat Z_\e^\pm$ is a supersolution with the Robin inequality above and $\widehat Z_\e^\pm(0)\ge|u_{\e,0}^\pm|$.
Applying \cref{pro:basic_estimates}(ii) to $Z_\e^\pm\mp u_\e^\pm$, where $Z_\e^\pm=\widehat Z_\e^\pm\circ\Theta_\e^{-1}$, gives
\[
 |u_\e^\pm|\le Z_\e^\pm\le(u^*+\tau F^*)\exp\bigl(K_\ell\tau+C(1+M)\ell\bigr)
 \quad\text{on }Q_\e^\pm(I_\tau).
\]
Choosing $T_{\mathrm{geom}}=\ell_0^2$ and $\ell=\sqrt\tau\le\ell_0$ gives $K_\ell\tau+C(1+M)\ell\le C\bigl((1+M)\sqrt\tau+(1+M)^2\tau\bigr)$, which is \eqref{eq:short_sup}.
\end{proof}

\begin{remark}
\Cref{pro:basic_estimates,lem:short_sup} apply to every weak solution of the coupled problem whose heights lie in $\mathcal H_{M}(I_\tau)$, since its concentrations coincide with the solutions of \cref{lem:prescribed_bulk} for that motion.
\end{remark}

With \cref{lem:short_sup}, we can ensure, after fixing $M$, $\|u_\e\|_\infty\le u^*+1$ for $\tau$ sufficiently small.
However, for the surface equation, the uniform $L^\infty$-bound alone is not sufficient; we also need spatial Hölder control of the concentration traces on the moving membranes.
Since the fast outer diffusion does not provide a uniform time modulus as $\e\to0$, we use a local boundary estimate that requires only boundedness of the coefficients and boundary data in time.
To that end, we now apply the local estimate of Lemma~\ref{lem:local_holder} to the cellwise pullbacks introduced in Section~\ref{ssec:rescaling}.
More precisely, after the rescaling \(x=\e(k+y)\), the restrictions of \(\widehat u_\e^\pm\) to the membrane collars are represented by \(\widehat\mu_{\e,k}^\pm\), and the fixed reference interface \(\Gamma\) allows us to derive uniform spatial Hölder bounds for their traces.

\begin{proposition}
\label{pro:moving_holder}
Let $\bh\in\mathcal H_M(I_\tau)$ with $\bh(0)=0$ and let $L>0$ such that the
corresponding solution satisfies $\|u_\e\|_\infty\le L$.
There exist $\alpha\in(0,1)$ and $C>0$, depending only on $M,L,U^*,F^*$ and the
reference geometry, such that for all $k\in\indexI$:
\begin{enumerate}[label=(\roman*)]
 \item \emph{Outer phase.}
 \begin{equation}\label{eq:mu_equation}
  \sup_{t\in I_\tau}\|\widehat\mu_{\e,k}^{+}(t)\|_{C^\alpha(\Gamma)}
  +\e^\alpha\|\widehat\mu_{\e,k}^{+}\|_{C^{\alpha/2}(I_\tau;C(\Gamma))}\le C.
 \end{equation}
 \item \emph{Inner phase.}
 \begin{equation}\label{eq:inner_holder}
  \|\widehat\mu_{\e,k}^{-}\|_{C^{\alpha/2,\alpha}(\overline{I_\tau}\times\overline{Y^{-}})}\le C.
 \end{equation}
\end{enumerate}
In particular,
\begin{equation}\label{eq:moving_holder}
 \max_{k,\pm}\sup_{t\in I_\tau}\|\widehat\mu_{\e,k}^{\pm}(t)\|_{C^\alpha(\Gamma)}\le C.
\end{equation}
\end{proposition}
\begin{proof}
Fix $k\in\indexI$ and suppress the cell index in the geometric coefficients.
\emph{Outer phase.}
On the cell scale, $\widehat\mu_{\e,k}^+$ diffuses with coefficient $\e^{-2}d^+$, see \eqref{eq:fixed_cell:plus}.
In the fast time $s=t/\e^2$, the function $z(s,y)=\widehat\mu_{\e,k}^{+}(\e^2s,y)$ solves, with $a_h=A_h/J_h$,
\begin{align*}
    \partial_s z-\dive_y(d^+a_h\nabla_y z)
    &=\bigl(\e^2b_h+d^+a_h\nabla_y\log J_h\bigr)\cdot\nabla_y z+\e^2\widehat f_{\e,k}^+
    &&\text{in }Y_\Sigma^+,\\
    d^+a_h\nabla_y z\cdot\nu^+&=\e^2m_h\beta(h)\partial_t h\,z
    &&\text{on }\Gamma,
\end{align*}
where $J_h^{-1}\dive_y(A_h\nabla_yz)=\dive_y(a_h\nabla_yz)+a_h\nabla_y\log J_h\cdot\nabla_yz$ was used.
By \eqref{ass:geometry_bootstrap} and \cref{lem:geometric_estimates}, $d^+a_h$ is uniformly elliptic and bounded, the drift is bounded by $C(1+M)$, the source by $F^*$, and the conormal datum by $CML$, uniformly in $\e$ and $k$.
By \ref{ass:initial_data}, the initial value $z_0=u_{\e,0}^+(\e(k+\cdot))$ satisfies $K_0\le CU^*$ in \cref{lem:local_holder}; since the negative-time extension there uses the source $-\Delta_yz_0$ in the variable $s$, no factor $\e^{-2}$ enters.
As $\Gamma$ has positive distance from $\Sigma$, there are fixed radii $0<\varrho<\varrho_0$ with $B_{\varrho_0}(x_0)\cap\partial Y_\Sigma^+=B_{\varrho_0}(x_0)\cap\Gamma$ for all $x_0\in\Gamma$.
\Cref{lem:local_holder} with $D=Y_\Sigma^+$ and $\Gamma_N=\Gamma$ thus gives $\|z\|_{C^{\alpha/2,\alpha}(\overline{Q_\varrho})}\le C$ for all $x_0\in\Gamma$ and $s_0\in(0,\e^{-2}\tau]$, with $\alpha$ and $C$ independent of $\e$, $k$, $x_0$, and $s_0$.
Since the cylinders have fixed size, these local bounds combine with $\|z\|_\infty\le L$ to a bound on all of $[0,\e^{-2}\tau]\times\Gamma$, as at the end of Step~2 in the proof of \cref{lem:local_holder}; the constant does not depend on the length of the interval.
Returning to $t=\e^2s$ converts the temporal Hölder bound into the factor $\e^{\alpha}$ in \eqref{eq:mu_equation}, while the spatial bound is unchanged.

\emph{Inner phase.}
Here no time rescaling is needed: $\widehat\mu_{\e,k}^-$ solves the same type of problem in $Y^-$ with diffusion $d^-a_h$, drift bounded by $C(1+M)$, source bounded by $F^*$, and conormal datum $-m_h\beta(h)\partial_t h\,\widehat\mu_{\e,k}^-$ bounded by $CML$; by \ref{ass:initial_data}, $K_0\le CU^*$.
The same argument with $D=Y^-$ gives the bound near $\Gamma$, and interior Hölder estimates \cite[Chapter~VI]{Lieberman1996} give it on the rest of $\overline{Y^-}$.
Taking the smallest of the exponents proves the assertion.
\end{proof}

\begin{remark}\label{rem:scaled_holder}
Since $\widehat\mu_{\e,k}^\pm(t,y)=\widehat u_\e^\pm(t,\e(k+y))$, \eqref{eq:moving_holder} is a Hölder bound at the cell scale: on a single physical interface $\Gamma_{\e,k}$, the unscaled seminorm $\|\widehat u_\e^\pm(t)\|_{C^\alpha(\Gamma_{\e})}$ may grow like $\e^{-\alpha}$.
\end{remark}

\subsection{Interface evolution with prescribed concentration forcing}
\label{ssec:schauder}
\Cref{pro:moving_holder} provides a uniform spatial Hölder bound for the concentration jump driving the membrane motion, whereas its temporal Hölder bound may still depend on $\e$ uncontrollably via
\[
    \e^\alpha\max_k
    \|\widehat\mu_{\e,k}^+\|_{
        C^{\alpha/2}(I_\tau;C(\Gamma))}
    \le C.
\]
In particular, we cannot infer any uniform Hölder bound in time for the concentration jumps $g_k=\widehat\mu_{\e,k}^--\widehat\mu_{\e,k}^+$.
Such a bound would allow us to apply classical parabolic Schauder theory directly to the interface evolution equation with uniform regularity estimates for the membrane heights.
To circumvent this missing control in time, we develop estimates that require only spatial Hölder regularity of the forcing.

Throughout this subsection we fix $0<\beta<\alpha<1$ and auxiliary parameters
\[
  \alpha':=\tfrac{\alpha+\beta}{2},\qquad s\in(\beta,\alpha'),\qquad
  q,m\in(1,\infty)\ \text{with}\ s-\tfrac2q-\tfrac2m>\beta .
\]
For $\tau\in(0, T^*)$, we set
\begin{align*}
  E_{\tau}&:=\bigl\{h\in W^{1,q}(I_\tau;W^{s,m}(\Gamma))\cap L^q(I_\tau;W^{s+2,m}(\Gamma)) : h(0)=0\bigr\},\\
  \mathcal G_{B,\alpha}(I_\tau)&:=\Bigl\{g\in C(\overline{I_\tau}\times\Gamma):
  \sup_{t\in\overline{I_\tau}}\|g(t)\|_{C^\alpha(\Gamma)}\le B\Bigr\}.
\end{align*}
We note that $\mathcal G_{B,\alpha}(I_\tau)$ encodes exactly the estimates that are available via \cref{pro:moving_holder}.
By the trace theorem for maximal-regularity spaces, $E_{\tau}\hookrightarrow C(\overline{I_\tau};C^{2+\beta}(\Gamma))$
with a constant $C_E$ independent of $\tau\le T^*$; see \cite[Theorem~6.1.11 and Section~6.1.5]{PruessSimonett2016}.
With the following lemma, we show that the spatial control alone is already sufficient to keep the membrane evolution uniformly bounded on a common short time interval.

\begin{lemma}
\label{lem:surface_spatial}
Let $B>0$.
There are $\varrho>0$ and $T_{\mathrm{surf}}=T_{\mathrm{surf}}(B)>0$ such that, for every $0<\tau\le T_{\mathrm{surf}}$ and $g\in\mathcal G_{B,\alpha}(I_\tau)$, the following holds:
\begin{enumerate}
  \item[$(i)$] The problem
  \begin{equation}\label{eq:surface_spatial}
    \beta(h)\partial_th+\sigma\widehat H[h]=\lambda g,\qquad h(0)=0,
  \end{equation}
  has a unique solution $h=S_{\tau}(g)$ in $\{h\in E_{\tau}:\|h\|_{E_{\tau}}\le\varrho\}$. It satisfies
  \[
    h\in\mathcal X^\beta(I_\tau),\qquad \|h\|_{C(\overline{I_\tau};C^{2+\beta}(\Gamma))}\le \eta^*/4,\qquad
    \partial_th(0)=-\sigma H_\Gamma+\lambda g(0).
  \]
  \item[$(ii)$] The solution operator $S_\tau\colon\mathcal G_{B,\alpha}(I_\tau)\to\mathcal X^\beta(I_\tau)$ is Lipschitz continuous when $\mathcal G_{B,\alpha}(I_\tau)$ is endowed with the norm of $C(\overline{I_\tau};C^{\alpha'}(\Gamma))$.
  \item[$(iii)$] $S_{\tau}$ maps subsets of $\mathcal G_{B,\alpha}(I_\tau)$ that are bounded in   $C^{\alpha/2,\alpha}(\overline{I_\tau}\times\Gamma)$ into relatively compact subsets of $\mathcal X^\beta(I_\tau)$.
\end{enumerate}
\end{lemma}

\begin{proof}
We use two standard facts on the closed surface $\Gamma$: for $s<s'<1$,
\begin{equation}\label{eq:surface_multiplier}
  \|\varphi v\|_{W^{s,m}(\Gamma)}\le C\|\varphi\|_{C^{s'}(\Gamma)}\|v\|_{W^{s,m}(\Gamma)},
  \qquad
  C^{\alpha'}(\Gamma)\hookrightarrow W^{s,m}(\Gamma),
\end{equation}
the latter since $s<\alpha'$. Moreover, $(\partial_t-\sigma\Delta_\Gamma)$ has maximal
$L^q$-regularity on $W^{s,m}(\Gamma)$, see \cite[Theorem~6.1.11, Section~6.1.5]{PruessSimonett2016};
we write $\mathcal L_{\tau}^{-1}$ for its solution operator with zero initial value. Since the
initial value vanishes,
\begin{equation}\label{eq:surface_mr_constants}
  \|\mathcal L_{\tau}^{-1}f\|_{E_{\tau}}\le C_{\mathcal L}\|f\|_{L^q(I_\tau;W^{s,m}(\Gamma))},
  \qquad
  \|h\|_{C(\overline{I_\tau};C^{2+\beta}(\Gamma))}\le C_E\|h\|_{E_{\tau}}\quad(h\in E_{\tau}),
\end{equation}
where $C_{\mathcal L},C_E$ are independent of $\tau\le T^*$.

In local charts, the mean curvature of a normal graph over $\Gamma$ is a quasilinear second-order operator in $h$, and $\beta(h)$ is a smooth function of $(h,\nabla_\Gamma h)$ with $\beta(0)=1$, see \cite[Section~2.5]{PruessSimonett2016}. 
By the choice of $\eta^*$ and \cref{lem:geometric_estimates}, we can therefore write, on $\{\|h\|_{C^2(\Gamma)}\le\eta^*/2\}$,
\[
  -\sigma\beta(h)^{-1}\widehat H[h]
  =\mathcal A(h,\nabla_\Gamma h):\nabla_\Gamma^2h+\mathcal F(h,\nabla_\Gamma h),
\]
with smooth $\mathcal A,\mathcal F$, where $\mathcal A$ is uniformly elliptic,
$\mathcal A(0,0):\nabla_\Gamma^2=\sigma\Delta_\Gamma$, and $\mathcal F(0,0)=-\sigma H_\Gamma$.
Hence \eqref{eq:surface_spatial} is equivalent to
\begin{equation}\label{eq:surface_reformulated}
  \partial_th-\sigma\Delta_\Gamma h=-\sigma H_\Gamma+\mathcal Q[h]+\lambda\beta(h)^{-1}g,
  \qquad h(0)=0,
\end{equation}
with the remainder
\[
  \mathcal Q[h]:=\bigl(\mathcal A(h,\nabla_\Gamma h)-\mathcal A(0,0)\bigr):\nabla_\Gamma^2h
  +\mathcal F(h,\nabla_\Gamma h)-\mathcal F(0,0).
\]
In particular, $\mathcal Q[0]=0$, and the coefficient of $\nabla_\Gamma^2h$ in $\mathcal Q$ vanishes at $h=0$.
Consequently, \eqref{eq:surface_reformulated} is a quasilinear perturbation of $(\partial_t-\sigma\Delta_\Gamma)$.

\medskip
\emph{Proof of (i).}
Fix $B>0$, $\tau\le T^*$, and $g\in\mathcal G_{B,\alpha}(I_\tau)$.
Since $g$ is continuous with values in $C(\Gamma)$ and bounded in $C^\alpha(\Gamma)$, interpolation gives $g\in C(\overline{I_\tau};C^{\alpha'}(\Gamma))$, so $g$ is strongly measurable with values in $W^{s,m}(\Gamma)$ by \eqref{eq:surface_multiplier}.
Choose $\varrho>0$ such that $C_E\varrho\le\eta^*/4$, and set
\[
  \mathbb B_\varrho
  :=\{h\in E_\tau:\|h\|_{E_\tau}\le\varrho\}.
\]
Then every $h\in\mathbb B_\varrho$ satisfies
\[
  \|h\|_{C(\overline{I_\tau};C^{2+\beta}(\Gamma))}
  \le C_E\varrho\le\eta^*/4\le\eta^*/2.
\]
Therefore, the geometric coefficients are well defined, and we define
\[
  \Phi_g\colon\mathbb B_\varrho\to E_\tau,\qquad
  \Phi_g(h):=\mathcal L_\tau^{-1}
  \bigl(-\sigma H_\Gamma+\mathcal Q[h]
        +\lambda\beta(h)^{-1}g\bigr).
\]

Let $h_1,h_2\in\mathbb B_\varrho$, $w:=h_1-h_2$ and write
\[
  \mathcal A_j:=\mathcal A(h_j,\nabla_\Gamma h_j),
  \qquad
  \mathcal F_j:=\mathcal F(h_j,\nabla_\Gamma h_j),
  \qquad
  \mathcal A_0:=\mathcal A(0,0).
\]
Then,
\begin{align}\label{eq:q_lipschitz}
  \mathcal Q[h_1]-\mathcal Q[h_2]
  =(\mathcal A_1-\mathcal A_0):\nabla_\Gamma^2w
   +(\mathcal A_1-\mathcal A_2):\nabla_\Gamma^2h_2
   +\mathcal F_1-\mathcal F_2.
\end{align}
Fix $s'\in(s,1)$. Smoothness of the coefficients and the
embedding in \eqref{eq:surface_mr_constants} give
\[
  \|\mathcal A_j-\mathcal A_0\|_{C(\overline{I_\tau};C^{s'})}
  \le C\varrho,
  \qquad
  \|\mathcal A_1-\mathcal A_2\|_{C(\overline{I_\tau};C^{s'})}
  \le C\|w\|_{E_\tau}.
\]
Thus, by \eqref{eq:surface_multiplier}, each of the first two terms on the right-hand side of \eqref{eq:q_lipschitz} has $L^q(I_\tau;W^{s,m}(\Gamma))$-norm at most $C\varrho\|w\|_{E_\tau}$.
For the last term, smoothness of $\mathcal F$ and the
chain rule, together with the uniform $C^2$ bounds on
$h_1,h_2$, yield
\[
  \|\mathcal F_1-\mathcal F_2\|_{C(\overline{I_\tau};C^1)}
  \le C\|w\|_{C(\overline{I_\tau};C^2)}
  \le C\|w\|_{E_\tau}.
\]
Using $C^1(\Gamma)\hookrightarrow W^{s,m}(\Gamma)$ and integrating in time, we then obtain
\begin{equation}\label{eq:Q_lipschitz}
  \|\mathcal Q[h_1]-\mathcal Q[h_2]\|_{L^q(I_\tau;W^{s,m})}
  \le C(\varrho+\tau^{1/q})\|w\|_{E_\tau}.
\end{equation}

For the forcing term, the same smoothness argument,
the product estimate in $C^\alpha(\Gamma)$, and
$C^\alpha(\Gamma)\hookrightarrow W^{s,m}(\Gamma)$ yield
\begin{equation}\label{eq:forcing_estimates}
\begin{aligned}
  \|\beta(h_1)^{-1}g\|_{L^q(I_\tau;W^{s,m})}
  &\le CB\tau^{1/q},\\
  \|(\beta(h_1)^{-1}-\beta(h_2)^{-1})g\|_{L^q(I_\tau;W^{s,m})}
  &\le CB\tau^{1/q}\|w\|_{E_\tau}.
\end{aligned}
\end{equation}
All constants above are independent of $\tau\le T^*$.

We now choose $\varrho$ and $\tau$ to obtain a contraction.
By maximal regularity and the preceding estimates,
\[
  \|\Phi_g(h_1)-\Phi_g(h_2)\|_{E_\tau}
  \le C_{\mathcal L}C
       \bigl(\varrho+(1+B)\tau^{1/q}\bigr)
       \|h_1-h_2\|_{E_\tau}.
\]
Moreover, $\mathcal Q[0]=0$ and $\beta(0)=1$ imply
\[
  \|\Phi_g(0)\|_{E_\tau}
  \le C_{\mathcal L}C(1+B)\tau^{1/q}.
\]
Choose $\varrho_0>0$ such that
$C_{\mathcal L}C\varrho_0\le1/4$, and set
\[
  \varrho:=\min\{\varrho_0,\eta^*/(4C_E)\}.
\]
Next choose $0<T_{\mathrm{surf}}(B)\le T^*$ such that
\[
  C_{\mathcal L}C(1+B)T_{\mathrm{surf}}^{1/q}
  \le\min\{1/4,\varrho/2\}.
\]
For $\tau\le T_{\mathrm{surf}}$, the map $\Phi_g$ has
Lipschitz constant at most $1/2$ on $\mathbb B_\varrho$,
and
\[
  \|\Phi_g(h)\|_{E_\tau}
  \le \|\Phi_g(h)-\Phi_g(0)\|_{E_\tau}
      +\|\Phi_g(0)\|_{E_\tau}
  \le \tfrac12\|h\|_{E_\tau}+\tfrac{\varrho}{2}
  \le\varrho.
\]
Banach's fixed-point theorem therefore gives a unique
fixed point $h=S_\tau(g)\in\mathbb B_\varrho$.
This is precisely the unique solution of \eqref{eq:surface_spatial} in $\mathbb B_\varrho$, and
\[
  \|h\|_{C(\overline{I_\tau};C^{2+\beta}(\Gamma))}
  \le C_E\varrho\le\eta^*/4.
\]

The right-hand side of \eqref{eq:surface_reformulated} and $\Delta_\Gamma h$ belong to $C(\overline{I_\tau}\times\Gamma)$, because $h\in C(\overline{I_\tau};C^{2+\beta}(\Gamma))$ and $g$ is continuous.
Hence $\partial_th\in C(\overline{I_\tau}\times\Gamma)$, so $h\in\mathcal X^\beta(I_\tau)$, and evaluating \eqref{eq:surface_reformulated} at $t=0$ with $h(0)=0$, $\mathcal Q[0]=0$, and $\beta(0)=1$ gives $\partial_th(0)=-\sigma H_\Gamma+\lambda g(0)$.

\medskip
\emph{Proof of (ii).}
Let $g_1,g_2\in\mathcal G_{B,\alpha}(I_\tau)$ and $h_j=S_{\tau}(g_j)$.
Then
\[
  h_1-h_2=\bigl(\Phi_{g_1}(h_1)-\Phi_{g_1}(h_2)\bigr)
  +\lambda\mathcal L_{\tau}^{-1}\bigl(\beta(h_2)^{-1}(g_1-g_2)\bigr).
\]
Using the contraction estimate for $\Phi_{g_1}$ and absorbing $\frac12\|h_1-h_2\|_{E_\tau}$ into the left-hand side, we obtain
\[
  \|h_1-h_2\|_{E_{\tau}}\le 2\lambda C_{\mathcal L}\|\beta(h_2)^{-1}(g_1-g_2)\|_{L^q(I_\tau;W^{s,m})}
  \le C\|g_1-g_2\|_{C(\overline{I_\tau};C^{\alpha'}(\Gamma))}.
\]
By \eqref{eq:surface_mr_constants} this controls $h_1-h_2$ in $C(\overline{I_\tau};C^{2+\beta}(\Gamma))$.
Subtracting the two equations \eqref{eq:surface_reformulated} and using the smoothness of $\mathcal A$, $\mathcal F$, and $\beta^{-1}$ also controls $\partial_t(h_1-h_2)$ in $C(\overline{I_\tau}\times\Gamma)$.
Thus $S_{\tau}$ is even Lipschitz from $C(\overline{I_\tau};C^{\alpha'}(\Gamma))$ into $\mathcal X^\beta(I_\tau)$.

\medskip
\emph{Proof of (iii).}
Let $\mathcal G'\subset\mathcal G_{B,\alpha}(I_\tau)$ be bounded in $C^{\alpha/2,\alpha}(\overline{I_\tau}\times\Gamma)$ by $K$, and let $h=S_{\tau}(g)$ with $g\in\mathcal G'$.
Then \eqref{eq:surface_spatial} reads in local charts as
\[
  \partial_th-\mathcal A(h,\nabla_\Gamma h):\nabla_\Gamma^2h=\mathcal F(h,\nabla_\Gamma h)+\lambda\beta(h)^{-1}g,
\]
which is uniformly parabolic.
By $(i)$, $h$ is bounded in $C(\overline{I_\tau};C^{2+\beta}(\Gamma))$ by $\eta^*/4$, and $\partial_th$ is bounded in $C(\overline{I_\tau}\times\Gamma)$.
Interpolating in space between these bounds shows that $h$ and $\nabla_\Gamma h$ are bounded in $C^{\beta_0/2,\beta_0}$ for some $\beta_0\in(0,\beta]$.
Hence the coefficient $\mathcal A(h,\nabla_\Gamma h)$ and the right-hand side are bounded in $C^{\beta_0/2,\beta_0}$, and the parabolic Schauder estimate on the closed surface $\Gamma$ with zero initial value bounds $h$ in $C^{1+\beta_0/2,2+\beta_0}$.
In particular, $h$ and $\nabla_\Gamma h$ are now bounded in $C^{\alpha'/2,\alpha'}$, so the coefficient and the right-hand side are bounded in $C^{\alpha'/2,\alpha'}$, recalling $\alpha'<\alpha$ and $g\in C^{\alpha/2,\alpha}$.
A second application of the Schauder estimate bounds $h$ in $C^{1+\alpha'/2,2+\alpha'}(\overline{I_\tau}\times\Gamma)$, uniformly over $\mathcal G'$.
Since $\alpha'>\beta$, the Arzel\`a--Ascoli theorem shows that this set is relatively compact in $\mathcal X^\beta(I_\tau)$.
\end{proof}

\begin{corollary}
\label{cor:prescribed_concentrations}
Let $B,L>0$.
Let $0<\tau\le T_{\mathrm{surf}}(B)$ and let $g_{\e,k}:=\widehat\mu^-_{\e,k}-\widehat\mu^+_{\e,k}\in\mathcal G_{B,\alpha}(I_\tau)$ with $\|\widehat\mu^\pm_{\e,k}\|_\infty\le L$ for all $k\in \indexI$.
Then, the interface evolution equations
\[
    \beta(h_{\e,k})\partial_t h_{\e,k}
    +\sigma \widehat H[h_{\e,k}]
    =
    \lambda\bigl(
        \widehat\mu_{\e,k}^-
        -\widehat\mu_{\e,k}^+
    \bigr),
    \qquad
    h_{\e,k}(0)=0,
\]
admit a unique family of solutions
$\bh_\e=(h_{\e,k})_{k\in\indexI}\in\mathcal H_{M_{\mathrm{out}}}(I_\tau)$ satisfying
\[
    \max_{k\in\indexI}
    \|h_{\e,k}\|_{C(\overline{I_\tau};C^{2+\beta}(\Gamma))}
    \le\frac{\eta^*}4,\qquad
    \max_k\|\partial_th_{\e,k}\|_{L^\infty}\le \frac{\sigma H^*+2\lambda L}{b_*}.
\]
For each fixed $\e>0$, the solution map is continuous with respect to convergence of the concentration traces in $C(\overline{I_\tau};C^{\alpha'}(\Gamma))$ and maps families with uniformly bounded $C^{\alpha/2,\alpha}$ traces into relatively compact subsets of $\mathcal X_\e^\beta(I_\tau)$.
\end{corollary}
\begin{proof}
\Cref{lem:surface_spatial} applies to each interface with a common existence time $T_{\mathrm{surf}}=T_{\mathrm{surf}}(B)$ which also gives the asserted $C^{2+\beta}$ bounds.
Since $\|h_{\e,k}\|_{C^2}\le\eta^*/4\le\eta^*/2$, the surface equation gives the estimate for $\partial_th_{\e,k}$.

The continuity and compactness assertions follow from \cref{lem:surface_spatial}, since $\indexI$ is finite for each fixed $\e$.
\end{proof}

\subsection{Proof of \cref{thm:prepared_schauder}}
\label{ssec:proof_wellposed}

We now combine the solution operators of \cref{lem:prescribed_bulk} and \cref{lem:surface_spatial} in a Schauder fixed-point argument.
Starting from a prescribed interface motion, we solve the bulk diffusion problem and use the resulting concentration jumps to determine a new interface
motion. 
We show that this procedure defines a continuous and compact self-map on a suitable closed convex set, with an existence interval independent of $\e$.

We fix the constants in the following order;
each depends only on those chosen before.
\begin{enumerate}
  \item $L:=u^*+1$ and $M:=b_*^{-1}(\sigma H^*+2\lambda L)$.
  \item $T_{\mathrm{bulk}}=T_{\mathrm{bulk}}(M)\le T_{\mathrm{geom}}$ by \cref{lem:short_sup},
        such that every $\bh\in\mathcal H_M$ with $\bh(0)=0$ implies $\|u_\e\|_\infty\le L$.
  \item $\alpha$ and $B=2C$ with $C$ as in \eqref{eq:moving_holder} by \cref{pro:moving_holder}, depending only on $M,L,U^*,F^*$.
  \item $\beta:=\alpha/2$.
  \item $T:=\min\{T^*,T_{\mathrm{bulk}},T_{\mathrm{surf}}(B)\}$ and $I:=I_T=(0,T)$.
\end{enumerate}
We also introduce the nonempty closed convex subset
\begin{equation*}
 K_{\e}=\bigl\{\bh\in  X_{\e}^\beta(I):\ h_k(0)=0,\quad
       \max_k\|h_k\|_{C(\overline I;C^{2+\beta}(\Gamma))}\le\frac{\eta^*}2,\quad
       \max_k\|\partial_th_k\|_{C(\overline I\times\Gamma)}\le M\bigr\}
\end{equation*}
and note that every $\bh\in K_{\e}$ satisfies \eqref{ass:geometry_bootstrap} and, therefore, $K_\e\subset\mathcal H_{M}(I)$.

Now, given $\bh\in K_{\e}$, we solve the global prescribed-motion bulk problem to get the corresponding $\widehat u^\pm_\e[\bh]$ and form $g_k[\bh]=(\widehat\mu_{\e,k}^{-}[\bh]-\widehat\mu_{\e,k}^{+}[\bh])$.
We already know that $g_k[\bh]\in\mathcal G_{B,\alpha}(I)$ via \cref{pro:moving_holder}.
Next, we solve \eqref{eq:surface_spatial} with $g_k[\bh]$ as data (via \cref{lem:surface_spatial}) and denote the resulting family by $\widetilde\bh=(\widetilde h_k)_k$.
Since the new heights have norm at most $\eta^*/4$, their velocities satisfy, directly from their equation,
\[
 \|\partial_t\widetilde h_k\|_\infty
 \le b_*^{-1}\bigl(\sigma H^*+2\lambda(u^*+1)\bigr)=M.
\]
Thus, this construction induces an operator
\[
\Phi_\e\colon K_{\e}\to K_{\e}\qquad \Phi_\e\colon\bh\mapsto u_\e^\pm[\bh]\mapsto(g_k[\bh])_k\mapsto \widetilde\bh 
\]

By \cref{pro:moving_holder}, the jumps satisfy
\begin{equation}\label{eq:jump_bounds_fixed_eps}
  g_k[\bh]\in\mathcal G_{B,\alpha}(I),
  \qquad
  \|g_k[\bh]\|_{C^{\alpha/2,\alpha}(\overline I\times\Gamma)}\le C_\e
  \qquad(\bh\in K_\e,\ k\in \indexI).
\end{equation}
Let $\bh_j\to\bh$ in $\mathcal X_\e^\beta(I)$.
Since $\mathcal X_\e^\beta(I)\hookrightarrow\mathcal X_\e(I)$, \cref{lem:app_bulk_continuity} gives convergence of the corresponding bulk solutions in $C(\overline I;L^2(\Omega_\e^\pm))\cap L^2(I;H^1(\Omega_\e^\pm))$ on both phases, and hence $g_k[\bh_j]\to g_k[\bh]$ in $L^2(I\times\Gamma)$ by continuity of the trace.
By \eqref{eq:jump_bounds_fixed_eps} and the Arzel\`a--Ascoli theorem, every subsequence of $(g_k[\bh_j])_j$ has a further subsequence converging uniformly on $\overline I\times\Gamma$, and its limit must be $g_k[\bh]$; thus the whole sequence converges uniformly.
The interpolation inequality
\[
  \|g\|_{C^{\alpha'}(\Gamma)}\le C\|g\|_{C(\Gamma)}^{1-\alpha'/\alpha}\|g\|_{C^\alpha(\Gamma)}^{\alpha'/\alpha},
\]
applied pointwise in time together with the $C^\alpha$ bound in $\mathcal G_{B,\alpha}(I)$, upgrades this to convergence in $C(\overline I;C^{\alpha'}(\Gamma))$.
Since $S_T$ is Lipschitz from $C(\overline I;C^{\alpha'}(\Gamma))$ into $\mathcal X^\beta(I)$ by \cref{lem:surface_spatial}(ii) and $I_\e$ is finite, $\Phi_\e(\bh_j)\to\Phi_\e(\bh)$ in $\mathcal X_\e^\beta(I)$.
Moreover, \eqref{eq:jump_bounds_fixed_eps} and \cref{lem:surface_spatial}(iii) show that $\{S_T(g_k[\bh]):\bh\in K_\e\}$ is relatively compact in $\mathcal X^\beta(I)$ for each $k$, so $\Phi_\e(K_\e)$ is relatively compact in $\mathcal X_\e^\beta(I)$.
As $K_\e$ is nonempty, closed, bounded, and convex, Schauder's fixed-point theorem yields a (not necessarily unique) fixed point $\bh\in K_\e$, on the interval $\overline I$.
In particular, $(\widehat u^\pm_\e,\bh)$ is a weak solution in the sense of \cref{def:weak_solution} on $I$; the bounds of \cref{lem:short_sup} hold by the remark following that lemma. It satisfies the $\e$-uniform bounds
\begin{equation*}
\|h_k\|_{C(\overline I;C^{2+\beta}(\Gamma))}\le\frac{\eta^*}2,\qquad 
\|\partial_th_k\|_\infty\le M,\qquad
\|u^\pm_\e\|_\infty\le L,\qquad
\sup_t\|\widehat\mu^\pm_{\e,k}(t)\|_{C^\alpha(\Gamma)}\le B.
\end{equation*}

This completes the proof of \cref{thm:prepared_schauder}.

\section{Homogenization}\label{sec:homogenization}

Throughout this section we impose \ref{ass:geometry}--\ref{ass:sources}, \ref{ass:macro_initial}, and \ref{ass:macro_sources}, and fix the weak solutions $(u_\e^\pm,\bh)$ on $I=(0,T)$ from \cref{thm:prepared_schauder}, with fixed exponents $\alpha$, $\beta$; we consider a sequence $E$ of scales $\e\downarrow0$, and subsequences are taken within $E$.
With the unfolding operators of Appendix~\ref{app:unfolding}, set $U_\e^\pm=\mathcal T_\e\widehat u_\e^\pm$ and $h_\e^{\rm unf}=\mathcal R_\e(h_{\e,k})_k$, write $D_\e=\Omega_\e^+$, and let $P_\e^+$ be the extension operator of \cref{lem:hom_extension}.

\paragraph{Strategy of the proof}
The uniform bounds of \cref{thm:prepared_schauder} give weak limits, but they do not prevent the cell geometries from oscillating from cell to cell, so the nonlinear quantities $J_h$, $A_h$, $\widehat H[h]$, $\theta^\pm(h)$, and $D(h)$ are not identified.
We circumvent that issue by regarding each inclusion as a cell problem driven by the outer concentration and proceed in four steps:
\Cref{ssec:hom_outer} proves strong compactness of the outer concentration and its membrane traces.
\Cref{ssec:hom_stability} proves an $L^2$-stability estimate for the driven cell problem that requires no time regularity of the exterior trace.
\Cref{ssec:hom_strong} combines both results with cellwise compactness to obtain strong convergence of the microstructure for a.e.\ $x$.
\Cref{ssec:hom_flux} identifies the effective flux and proves \cref{thm:homogenization}.

\subsection{Outer compactness and preliminary cell limits}
\label{ssec:hom_outer}
The extended outer initial data converge strongly to $u_0^+$ in $L^2(\Omega)$ by \ref{ass:macro_initial}, and $f_\e^+\to f^+$ strongly in $L^2(I^*\times\Omega)$ by \ref{ass:macro_sources}.
For the inner phase, the source enters the reference problem of the $k$th cell as $f_{\e,k}^-\circ\Theta_{h_{\e,k}}$, where $f_{\e,k}^-(t,y)=f_\e^-(t,\e(k+y))$ is defined on the whole cell $Y$.
The corresponding unfolded sources satisfy
\begin{equation}\label{eq:source_unfold_limit}
 F_\e^{-}:=\mathcal T_\e f_\e^{-}
 \longrightarrow f^{-}
 \quad\text{in }L^p\bigl(I\times\Omega;C(\overline{Y})\bigr)
\end{equation}
for every finite $p$, by \ref{ass:macro_sources}, the uniform cellwise Lipschitz bound in \ref{ass:sources}, and interpolation; see \cref{rem:assumptions}.

The extension operators and the uniform time-translate estimate for the Jacobians established in
\cref{lem:hom_extension,lem:jac_time_translates}, together with the uniform bounds of
\cref{thm:prepared_schauder}, yield the following compactness result for the outer concentration.

\begin{proposition}
\label{thm:outer_compactness}
There exist a subsequence and a limit
\[
 u^{+}\in L^2(I;H^1(\Omega))\cap L^\infty(I\times\Omega)
\]
such that
\begin{align}
 P_\e^+\widehat u_\e^{+}&\to u^{+} &&\text{strongly in }L^p(I\times\Omega),\label{eq:outer_strong}\\
P_\e^+\widehat u_\e^{+}&\rightharpoonup u^{+} &&\text{weakly in }L^2(I;H^1(\Omega)),\notag\\
U_\e^{+}&\to u^{+}
 &&\text{strongly in }L^2(I\times\Omega\times\Gamma)\label{eq:outer_trace_strong}
\end{align}
for every $p\in[1,\infty)$.
\end{proposition}
\begin{proof}
Since the conservative equation controls the time derivative of the weighted mass rather than the concentration itself, we combine this control with the temporal modulus of the Jacobian to obtain time translates of the concentration, and then use spatial trace regularity to pass to the membrane values.
The energy estimate in \cref{pro:basic_estimates}(i), transferred to the reference geometry using \cref{lem:geometric_estimates} and extended by \cref{lem:hom_extension}, gives $\|v_\e\|_{L^2_tH^1_x}\le C$.
The $L^\infty$ bound follows from \cref{thm:prepared_schauder} and the $L^\infty$-bounded extension in \cref{lem:hom_extension}.
Set, extending by zero outside $D_\e$,
\[
 a_\e=\chi_{D_\e}J_\e,\qquad m_\e=a_\e v_\e,\qquad
 Q_\e=\chi_{D_\e}
       (d^{+}A_\e\nabla\widehat u_\e^{+}+J_\e b_\e\widehat u_\e^{+}).
\]
The conservative equation \eqref{eq:fixed_conservative} and zero total normal flux imply
\[
 \partial_t m_\e=\dive Q_\e+
       \chi_{D_\e}J_\e\widehat f_\e^{+}
       \quad\hbox{in }H^1(\Omega)',\qquad
 \|\partial_t m_\e\|_{L^2_tH^1_x{}'}\le C.
\]
Here $b_\e=O(\e)$ and no bound on $\partial_tJ_\e$ is used.
Writing $\Delta_s z(t)=z(t+s)-z(t)$, the algebraic identity
\[
 a_\e(t+s)|\Delta_s v_\e|^2
 =\Delta_s m_\e\,\Delta_s v_\e
       -\Delta_s a_\e\,v_\e(t)\Delta_s v_\e
\]
can be integrated in space and time.
The first term is bounded by $Cs$ using the duality between $H^1(\Omega)'$ and $H^1(\Omega)$; the second is bounded by $Cs^{1/2}$ by \eqref{eq:jac_time_modulus}.
Since $J_\e\ge c>0$ on $D_\e$ and the extension is linear, time independent, and $L^2$ bounded,
\begin{equation}\label{eq:outer_time_translate}
 \|\Delta_s v_\e\|_{L^2((0,T-s)\times\Omega)}^2
       \le C(s+s^{1/2}).
\end{equation}
The compactness criterion for time translates and the compact embedding $H^1(\Omega)\Subset L^2(\Omega)$ now give strong $L^2$ convergence \cite{Simon1987}.
Interpolation with the uniform $L^\infty$ bound gives \eqref{eq:outer_strong}.

By \cref{lem:unfold_convergence}, unfolding a strongly convergent $L^2$ sequence on $\Omega$ gives strong convergence to the same function, independent of $y$.
Moreover,
\[
 \|\nabla_yU_\e^{+}\|_{L^2(I\times\Omega\times Y^{+})}
       \le C\e\|\nabla\widehat u_\e^{+}\|_{L^2(I\times D_\e)}
       \longrightarrow0.
\]
Continuity of the trace operator gives strong $L^2$ convergence on $\Gamma$.
Applying this estimate pointwise in $(t,x)$ and then integrating proves \eqref{eq:outer_trace_strong}.
\end{proof}

We first record compactness of the individual reference-cell trajectories, which will strengthen the preliminary weak limits below.

\begin{lemma}
\label{lem:cell_paths}
Let $\beta'\in(0,\beta)$, where $\beta$ is the exponent in \cref{thm:prepared_schauder}.
\begin{enumerate}[label=(\roman*)]
 \item The set $\{\widehat\mu_{\e,k}^{-}:\e\in E,\ k\in\indexI\}$ is relatively compact in $C(\overline I\times\overline{Y^{-}})$ and weakly relatively compact in $L^2(I;H^1(Y^{-}))$.
 \item The set $\{h_{\e,k}:\e\in E,\ k\in\indexI\}$ is relatively compact in
 $C(\overline I;C^{2+\beta'}(\Gamma))$.
\end{enumerate}
\end{lemma}

\begin{proof}
$(i)$ By \cref{pro:moving_holder}, the functions $\widehat\mu_{\e,k}^{-}$ are bounded in $C^{\alpha/2,\alpha}(\overline I\times\overline{Y^{-}})$, uniformly in $\e$ and $k$.
The Arzel\`a--Ascoli theorem gives the first assertion.
The energy estimate of \cref{pro:basic_estimates}(i), applied to a single inclusion, has a right-hand side of order $\e^3$; after rescaling it bounds $\widehat\mu_{\e,k}^{-}$ in $L^2(I;H^1(Y^{-}))$, uniformly in $\e$ and $k$, which gives weak relative compactness.

$(ii)$ By \cref{thm:prepared_schauder}, the heights are bounded in $L^\infty(I;C^{2+\beta}(\Gamma))$ and in $W^{1,\infty}(I;C(\Gamma))$, uniformly in $\e$ and $k$. By interpolation they are equicontinuous with values in $C^{2+\beta'}(\Gamma)$,
and the compact embedding $C^{2+\beta}(\Gamma)\hookrightarrow C^{2+\beta'}(\Gamma)$ together
with the Arzel\`a--Ascoli theorem establishes compactness.
\end{proof}

\subsection{Stability of the driven cell problem}
\label{ssec:hom_stability}
For given initial profile $a\in W^{1,\infty}(Y^-)$, exterior trace $g\in L^\infty(I\times\Gamma)$, and source $F\in L^\infty(I;W^{1,\infty}(Y))$ prescribed on the whole cell, we consider the cell problem
\begin{empheq}[
    left={\systemlabel{sys:driven_cell}{\widehat P_{\rm cell}}\ \empheqlbrace}
]{alignat=2}
    \partial_t(J_h\widehat u^-)
    -\dive_y\bigl(d^{-}A_h\nabla_y\widehat u^-+J_hb_h\widehat u^-\bigr)
    &=J_h\,(F\circ\Theta_h)
    &&\text{in }I\times Y^{-},
    \label{eq:driven_cell_bulk}\\
    \bigl(d^{-}A_h\nabla_y\widehat u^-+J_hb_h\widehat u^-\bigr)\cdot n_\Gamma
    &=0
    &&\text{on }I\times\Gamma,
    \label{eq:driven_cell_flux}\\
    \beta(h)\partial_th+\sigma \widehat H[h]
    &=\lambda(\widehat u^--g)
    &&\text{on }I\times\Gamma,
    \label{eq:driven_cell_surface}\\
    \widehat u^-(0)&=a,\quad h(0)=0,
    \label{eq:driven_cell_initial}
\end{empheq}
which is the transformed version of a one-phase moving-boundary diffusion problem driven by an imposed exterior concentration $g$.
The source is evaluated in the moving inclusion $Y^-[h]\subset Y$, which is why it enters through $F\circ\Theta_h$.
Fix $M,L>0$.
We call $(\widehat u^-,h)$ a solution if
\[
 \widehat u^-\in L^2(I;H^1(Y^-))\cap L^\infty(I\times Y^-),\qquad
 h\in C(\overline{I};C^2(\Gamma))\cap W^{1,\infty}(I;C(\Gamma)),
\]
solves
\begin{multline}\label{eq:driven_weak}
 -\int_0^T\!\int_{Y^{-}}J_h\widehat u^-\,\partial_t\varphi
 +\int_0^T\!\int_{Y^{-}}(d^{-}A_h\nabla \widehat u^-+J_hb_h\widehat u^-)\cdot\nabla\varphi\\
 =\int_0^T\!\int_{Y^{-}}J_h\,(F\circ\Theta_h)\,\varphi+\int_{Y^{-}}a\varphi(0)
\end{multline}
for smooth tests vanishing at $T$ and \eqref{eq:driven_cell_surface} holds almost everywhere together with $h(0)=0$.
We say  that $(\widehat u^-,h)$ is an admissible solution with bounds $(M,L)$ if it additionally satisfies the estimates
\[
 \|\widehat u^-\|_{L^\infty(I\times Y^-)}\le L,\qquad
 \|h\|_{C(\overline{I};C^2(\Gamma))}\le\frac{\eta^*}2,\qquad
 \|\partial_th\|_{L^\infty(I\times\Gamma)}\le M.
\]

It is useful to compare the reference mass densities rather than the concentration pullbacks.
We set
\begin{equation}\label{eq:density_coefficients}
 q=J_h\widehat u^-,\qquad K_h=d^{-}A_h/J_h,\qquad
 c_h=b_h-K_h\nabla\log J_h
\end{equation}
and see that
\begin{equation}\label{eq:density_problem}
 \partial_tq-\dive(K_h\nabla q+c_hq)=J_h\,(F\circ\Theta_h),\qquad
 (K_h\nabla q+c_hq)\cdot n_\Gamma=0,\qquad q(0)=a.
\end{equation}

\begin{lemma}
\label{lem:density_integrability}
Let $(\widehat u^-,h)$ be an admissible solution of \systemref{sys:driven_cell}{\widehat P_{\rm cell}}.
Then, for every finite $p\ge2$,
\begin{equation}\label{eq:density_lp}
 \|\nabla(J_h\widehat u^-)\|_{L^p(I\times Y^{-})}\le C_p(M,L,\|a\|_{W^{1,\infty}},\|F\|_\infty).
\end{equation}
\end{lemma}
\begin{proof}
We want to apply the conormal estimate from \cite[Theorem~2.2]{DongZhang2016} for scalar equations with coefficients that are merely measurable in time.
We write $Q_T=I\times Y^{-}$ and fix $p\ge2$.

By admissibility and \cref{lem:geometric_estimates}, we have that $K_h$ is uniformly elliptic with constant $\kappa_0$, $\|K_h\|_{L^\infty(I;W^{1,\infty}(Y^-))}\le L_K$, $\|c_h\|_\infty\le B_c$, and $\|q\|_\infty\le CL$, with constants depending only on $M$, $L$, $d^-$, and the reference geometry.
Since $K_h(t,\cdot)$ is Lipschitz uniformly in $t$, its mean oscillation on balls of radius $R$ is at most $CL_KR$, so the small-oscillation condition \cite[Assumption~2.1]{DongZhang2016} holds on balls of radius $R_0=R_0(p,\kappa_0,L_K)$.

The function $v=q-a$ satisfies $v(0)=0$ and, weakly,
\begin{gather*}
 \partial_tv-\dive(K_h\nabla v)=\dive\mathcal G+\mathcal F,\qquad
 (K_h\nabla v+\mathcal G)\cdot n_\Gamma=0,\\
 \mathcal G=K_h\nabla a+c_hq,\qquad \mathcal F=J_h(F\circ\Theta_h),
\end{gather*}
with $\|\mathcal G\|_{L^p(Q_T)}+\|\mathcal F\|_{L^p(Q_T)}\le C(\|\nabla a\|_\infty+\|q\|_\infty+\|F\|_\infty)$.
Since $v(0)=0$ in $L^2(Y^-)$, extending $v$, $\mathcal G$, and $\mathcal F$ by zero and $K_h$ by $d^-I$ to $t<0$ gives a weak solution on $(-\infty,T)$ in the sense of \cite{DongZhang2016}.

For $\Lambda$ larger than the thresholds of \cite[Theorem~2.2]{DongZhang2016} for $s=2$ and $s=p$, the function $w=e^{-\Lambda t}v$ solves $\partial_tw-\dive(K_h\nabla w)+\Lambda w=\dive(e^{-\Lambda t}\mathcal G)+e^{-\Lambda t}\mathcal F$.
The theorem gives a solution $w_p$ with $\nabla w_p\in L^p$, which vanishes for $t<0$; as $p\ge2$, it lies in the same $L^2$ energy class as $w$, so uniqueness for $s=2$ gives $w=w_p$.
Hence
\[
 \|\nabla q\|_{L^p(Q_T)}\le\|\nabla a\|_{L^p(Q_T)}+C_pe^{\Lambda T}\bigl(\|\mathcal G\|_{L^p(Q_T)}+\|\mathcal F\|_{L^p(Q_T)}\bigr),
\]
which proves \eqref{eq:density_lp}.
\end{proof}

\begin{theorem}[Stability of the driven cell problem]
\label{thm:cell_stability}
Let $(\widehat u_j^-,h_j)$, $j=1,2$, be admissible solutions with common bounds $(M,L)$ for the inputs 
\[
(a_j,g_j,F_j)\in W^{1,\infty}(Y^-)\times L^\infty(I\times\Gamma)\times L^\infty(I;W^{1,\infty}(Y)).
\]
Then,
\begin{multline}
 \sup_{0\le t\le T}\bigl(\|Q(t)\|_{L^2(Y^{-})}^2
                          +\|\rho(t)\|_{H^1(\Gamma)}^2\bigr)
   +\int_0^T\bigl(\|\nabla Q\|_{L^2(Y^{-})}^2
              +\|\partial_t\rho\|_{L^2(\Gamma)}^2
              +\|\rho\|_{H^2(\Gamma)}^2\bigr)\di t\notag\\
 \le C\left(\|a_1-a_2\|_{L^2(Y^{-})}^2
       +\|g_1-g_2\|_{L^2(I\times\Gamma)}^2
       +\|F_1-F_2\|_{L^2(I\times Y)}^2\right),
 \label{eq:cell_stability}
\end{multline}
where
\[
 Q=J_{h_1}\widehat u^-_1-J_{h_2}\widehat u^-_2,\qquad \rho=h_1-h_2.
\]
The constant depends only on $M,L$, the bounds on the inputs $(a_j,g_j,F_j)$, $T^*$, the fixed parameters $d^-,\sigma,\lambda$, and the reference geometry.
In particular, the driven cell solution is unique in this class.
\end{theorem}
\begin{proof}
We compare the reference mass densities to avoid differentiating the Jacobians in time and combine their energy estimate with the surface estimates.
For convenience, we write $J_j=J_{h_j}$, $K_j=K_{h_j}$, $c_j=c_{h_j}$, $q_j=J_j\widehat u^-_j$, $\delta g=g_1-g_2$, and $\delta F=F_1-F_2$.
By \cref{lem:app_coefficient_differences}, we have the Lipschitz estimates
\begin{align}
 \|K_1-K_2\|_{L^3(Y^{-})}
       &\le C\|\rho\|_{W^{1,3}(\Gamma)},\\
 \|c_1-c_2\|_{L^2(Y^{-})}
       &\le C\bigl(\|\rho\|_{H^2(\Gamma)}
                      +\|\partial_t\rho\|_{L^2(\Gamma)}\bigr),
       \label{eq:cell_coeff_diff}\\
 \|J_1-J_2\|_{L^2(Y^{-})}+
 \|J_1-J_2\|_{L^2(\Gamma)}
 +\|\beta(h_1)-\beta(h_2)\|_{L^2(\Gamma)}
       &\le C\|\rho\|_{H^1(\Gamma)}.\label{eq:stab_coeff}
\end{align}

\emph{Step 1: Surface estimate.}
Dividing both membrane laws by $\beta(h_j)$ and using the quasilinear form from the proof of \cref{lem:surface_spatial}, with $\mathcal A_j:=\mathcal A(h_j,\nabla_\Gamma h_j)$ and $\mathcal F_j:=\mathcal F(h_j,\nabla_\Gamma h_j)$, we obtain
\[
  \partial_t\rho-\mathcal A_1:\nabla^2_\Gamma\rho=R,\qquad\rho(0)=0,
\]
where
\begin{multline*}
  R:=(\mathcal A_1-\mathcal A_2):\nabla^2_\Gamma h_2+\mathcal F_1-\mathcal F_2
  +\lambda\bigl(\beta(h_1)^{-1}-\beta(h_2)^{-1}\bigr)\bigl(\widehat u_2^--g_2\bigr)\\
  +\lambda\beta(h_1)^{-1}\Bigl(\frac{Q}{J_1}-\frac{q_2(J_1-J_2)}{J_1J_2}\Bigr)
  -\lambda\beta(h_1)^{-1}\delta g.
\end{multline*}
The bounds on $h_1,h_2$ in $C^2$ and on $\widehat u_2^-$, $g_2$ in $L^\infty$, together with \eqref{eq:stab_coeff}, give
\[
  \|R(t)\|_{L^2(\Gamma)}\le C\bigl(\|\rho(t)\|_{H^1(\Gamma)}+\|Q(t)\|_{L^2(\Gamma)}+\|\delta g(t)\|_{L^2(\Gamma)}\bigr).
\]
Since $\mathcal A_1$ is symmetric, uniformly elliptic, and bounded in
$L^\infty(I;W^{1,\infty}(\Gamma))$, we have
\[
  \int_\Gamma
    (\mathcal A_1:\nabla_\Gamma^2\rho)\,\Delta_\Gamma\rho
  \ge c\|\nabla_\Gamma^2\rho\|_{L^2(\Gamma)}^2
      -C\|\rho\|_{H^1(\Gamma)}^2.
\]
Indeed, localization in smooth charts and integration by parts yield a coercive principal term, while derivatives of $\mathcal A_1$, commutators, and curvature contributions give lower-order terms controlled by
\[
  C\|\rho\|_{H^1(\Gamma)}\|\rho\|_{H^2(\Gamma)}.
\]
For the corresponding localization and commutator arguments for elliptic operators on closed hypersurfaces, we refer to \cite[Section~6.4]{PruessSimonett2016}.
The lower-order terms are absorbed by Young's inequality.
Testing the surface equation with $\rho-\Delta_\Gamma\rho$, using this coercivity estimate and Young's inequality, we obtain
\[
  \frac{\di}{\di t}\|\rho\|_{H^1(\Gamma)}^2
  +c\|\rho\|_{H^2(\Gamma)}^2
  \le C\bigl(
    \|\rho\|_{H^1(\Gamma)}^2
    +\|R\|_{L^2(\Gamma)}^2
  \bigr).
\]
Moreover, the equation itself gives
\[
  \|\partial_t\rho\|_{L^2(\Gamma)}^2
  \le C\bigl(
    \|\rho\|_{H^2(\Gamma)}^2
    +\|R\|_{L^2(\Gamma)}^2
  \bigr).
\]
Integrating in time, recalling $\rho(0)=0$, and using the preceding bound on $R$, we conclude that
\begin{equation*}
  \|\rho(t)\|_{H^1(\Gamma)}^2
  +\int_0^t
    \bigl(
      \|\partial_s\rho\|_{L^2(\Gamma)}^2
      +\|\rho\|_{H^2(\Gamma)}^2
    \bigr)\di s
  \le C\int_0^t
    \bigl(
      \|Q\|_{L^2(\Gamma)}^2
      +\|\rho\|_{H^1(\Gamma)}^2
      +\|\delta g\|_{L^2(\Gamma)}^2
    \bigr)\di s.
\end{equation*}
Combined with the trace inequality $\|Q\|_{L^2(\Gamma)}^2\le\vartheta\|\nabla Q\|_{L^2(Y^-)}^2+C_\vartheta\|Q\|_{L^2(Y^-)}^2$, this gives, for every $\vartheta>0$,
\begin{multline}\label{eq:surface_difference_energy}
  \|\rho(t)\|_{H^1(\Gamma)}^2+\int_0^t\bigl(\|\partial_s\rho\|_{L^2(\Gamma)}^2+\|\rho\|_{H^2(\Gamma)}^2\bigr)\di s\\
  \le C_1\int_0^t\bigl(\vartheta\|\nabla Q\|_{L^2(Y^-)}^2+C_\vartheta\|Q\|_{L^2(Y^-)}^2
  +\|\rho\|_{H^1(\Gamma)}^2+\|\delta g\|_{L^2(\Gamma)}^2\bigr)\di s.
\end{multline}

\medskip
\emph{Step 2: Bulk estimate.}
We subtract \eqref{eq:density_problem} for $q_1,q_2$ and get
\begin{multline*}
  \partial_tQ-\dive\bigl(K_1\nabla Q+c_1Q+(K_1-K_2)\nabla q_2+(c_1-c_2)q_2\bigr)\\
  =J_1(\delta F\circ\Theta_{h_1})
  +J_1\bigl(F_2\circ\Theta_{h_1}-F_2\circ\Theta_{h_2}\bigr)
  +(J_1-J_2)(F_2\circ\Theta_{h_2})
\end{multline*}
with zero total conormal flux on $\Gamma$ and $Q(0)=a_1-a_2$.
Since $Q\in L^2(I;H^1(Y^-))$ and $\partial_tQ\in L^2(I;H^1(Y^-)')$, we may test the weak form with $Q$, see \cite[Chapter~XVIII, §1]{DautrayLions1992}.
A standard energy estimate using the ellipticity of $K_1$ and Young's inequality gives
\begin{multline*}
  \frac{\mathrm d}{\mathrm dt}\|Q\|_{L^2(Y^-)}^2+\kappa_0\|\nabla Q\|_{L^2(Y^-)}^2
  \le C\Bigl(\|Q\|_{L^2(Y^-)}^2+\|\rho\|_{H^2(\Gamma)}^2+\|\partial_t\rho\|_{L^2(\Gamma)}^2\\
  +\|\delta F\|_{L^2(Y)}^2+\|(K_1-K_2)\nabla q_2\|_{L^2(Y^-)}^2\Bigr),
\end{multline*}
where we used \eqref{eq:stab_coeff}, the boundedness of $J_1$, $q_2$, and $F_2$, the change of variables $y\mapsto\Theta_{h_1}(y)$, and the Lipschitz bound
\[
  \|F_2\circ\Theta_{h_1}-F_2\circ\Theta_{h_2}\|_{L^2(Y^-)}
  \le\|\nabla_yF_2\|_{L^\infty}\|\Theta_{h_1}-\Theta_{h_2}\|_{L^2(Y^-)}
  \le C\|\rho\|_{L^2(\Gamma)},
\]
which follows from $\Theta_{h_1}-\Theta_{h_2}=\mathcal E\rho$ (see Appendix~\ref{app:hanzawa}).
Set $G(t):=\|\nabla q_2(t)\|_{L^6(Y^-)}$.
Hölder's inequality, \eqref{eq:stab_coeff}, the Gagliardo--Nirenberg inequality on $\Gamma$, and Young's inequality give
\begin{align*}
  \|(K_1-K_2)\nabla q_2\|_{L^2(Y^-)}^2
  &\le CG^2\|\rho\|_{W^{1,3}(\Gamma)}^2
  \le CG^2\|\rho\|_{H^1(\Gamma)}^{4/3}\|\rho\|_{H^2(\Gamma)}^{2/3}\\
  &\le C\bigl(\|\rho\|_{H^2(\Gamma)}^2+G^3\|\rho\|_{H^1(\Gamma)}^2\bigr),
\end{align*}
and $G^3\in L^1(I)$ by \cref{lem:density_integrability} for $p=6$.
Integrating in time, we obtain constants $C_0,C>0$ with
\begin{multline}\label{eq:stab_bulk}
  \|Q(t)\|_{L^2(Y^-)}^2+\kappa_0\int_0^t\|\nabla Q\|_{L^2(Y^-)}^2\di s\\
  \le\|a_1-a_2\|_{L^2(Y^-)}^2
  +C_0\int_0^t\bigl(\|\partial_s\rho\|_{L^2(\Gamma)}^2+\|\rho\|_{H^2(\Gamma)}^2\bigr)\di s\\
  +C\int_0^t\bigl(\|Q\|_{L^2(Y^-)}^2+(1+G^3)\|\rho\|_{H^1(\Gamma)}^2+\|\delta F\|_{L^2(Y)}^2\bigr)\di s.
\end{multline}

\medskip
\emph{Step 3: Closing the estimate.}
We multiply \eqref{eq:surface_difference_energy} by $C_0+1$, add \eqref{eq:stab_bulk}, and choose $\vartheta>0$ such that $(C_0+1)C_1\vartheta\le\kappa_0/2$.
Then the gradient term on the right is absorbed, and so is the $C_0$-term from \eqref{eq:stab_bulk}.
For $E(t):=\|Q(t)\|_{L^2(Y^-)}^2+\|\rho(t)\|_{H^1(\Gamma)}^2$, this yields
\begin{equation*}
  E(t)+\frac{\kappa_0}{2}\int_0^t\|\nabla Q\|_{L^2(Y^-)}^2\di s+\int_0^t\bigl(\|\partial_s\rho\|_{L^2(\Gamma)}^2+\|\rho\|_{H^2(\Gamma)}^2\bigr)\di s
  \le CD+C\int_0^t(1+G^3)E\di s,
\end{equation*}
where 
\[
D:=\|a_1-a_2\|_{L^2(Y^{-})}^2
       +\|g_1-g_2\|_{L^2(I\times\Gamma)}^2
       +\|F_1-F_2\|_{L^2(I\times Y)}^2.
\]
Since $1+G^3\in L^1(I)$, Gronwall's inequality gives $E(t)\le CD\exp\bigl(C\int_0^T(1+G^3)\bigr)\le CD$, and inserting this into the right-hand side bounds the integral terms as well.
Equal inputs give $D=0$, hence $Q=\rho=0$, and $\widehat u_1^-=\widehat u_2^-$ since $J_{h_1}=J_{h_2}>0$.
\end{proof}

\subsection{Strong convergence of the microstructure}
\label{ssec:hom_strong}
We now combine the outer trace convergence with the stability of the driven cell problem to prove the inner and geometric convergence assertions of \cref{thm:homogenization}.
\begin{lemma}\label{lem:cell_closed}
Let $(\widehat u_n^-,h_n)$ be admissible solutions of
\systemref{sys:driven_cell}{\widehat P_{\rm cell}} with common bounds $(M,L)$ for inputs $(a_n,g_n,F_n)$.
Assume that $(\nabla\widehat u_n^-)_n$ is bounded in $L^2(I\times Y^-)$ and
\begin{align*}
 a_n&\to a \ \text{ in }C(\overline{Y^-}), &
 g_n&\to g \ \text{ in }L^2(I\times\Gamma), &
 F_n&\to F \ \text{ in }L^2(I\times Y),\\
 \widehat u_n^-&\to\widehat u^- \ \text{ in }C(\overline I\times\overline{Y^-}), &
 h_n&\to h \ \text{ in }C(\overline I;C^{2+\beta'}(\Gamma)). &&
\end{align*}
Then $(\widehat u^-,h)$ is an admissible solution of
\systemref{sys:driven_cell}{\widehat P_{\rm cell}} with bounds $(M,L)$ for the inputs
$(a,g,F)$, which satisfy the same bounds as $(a_n,g_n,F_n)$.
Moreover, $\nabla\widehat u_n^-\rightharpoonup\nabla\widehat u^-$ in $L^2(I\times Y^-)$ and
$\partial_th_n\to\partial_th$ in $L^2(I\times\Gamma)$.
\end{lemma}
\begin{proof}
By weak-$*$ compactness, the limits $(a,g,F)$ satisfy the same bounds as
$(a_n,g_n,F_n)$, and the bounds $\|\widehat u^-\|_\infty\le L$ and
$\|h\|_{C(\overline I;C^2(\Gamma))}\le\eta^*/2$ pass to the limit.
Since $(\nabla\widehat u_n^-)_n$ is bounded in $L^2(I\times Y^-)$, every
subsequence has a weakly convergent subsequence, and since
$\widehat u_n^-\to\widehat u^-$ uniformly, its limit is $\nabla\widehat u^-$.
Since $\widehat H[h_n]$ and $\beta(h_n)^{-1}$ also converge uniformly,
\eqref{eq:driven_cell_surface} gives
\[
 \partial_th_n=\beta(h_n)^{-1}\bigl[-\sigma\widehat H[h_n]+\lambda(\widehat u_n^--g_n)\bigr]
 \to\beta(h)^{-1}\bigl[-\sigma\widehat H[h]+\lambda(\widehat u^--g)\bigr]
 \quad\text{in }L^2(I\times\Gamma).
\]
The limit is the distributional time derivative of $h$, so $h$ satisfies
\eqref{eq:driven_cell_surface} and $h(0)=0$; by weak-$*$ compactness,
$\|\partial_th\|_{L^\infty}\le M$.
The coefficients $J_{h_n}$, $A_{h_n}$ converge uniformly, and $b_{h_n}\to b_h$ in $L^2$
As $F_n$ is uniformly Lipschitz in $y$ and $\Theta_{h_n}\to\Theta_h$ uniformly, $F_n\circ\Theta_{h_n}\to F\circ\Theta_h$ in $L^2(I\times Y^-)$.
Together with the weak convergence of $\nabla\widehat u_n^-$, this passes every term of
\eqref{eq:driven_weak} to the limit.
\end{proof}

\begin{proposition}
\label{thm:micro_compactness}
There are limits $\widehat u^{-}\in L^2(I\times\Omega;H^1(Y^-))$ and $h\in L^p(\Omega;C(\overline{I};C^{2+\beta'}(\Gamma)))$ such that, for every finite $p$ and $0<\beta'<\beta$,
\begin{align}
 U_\e^{-}&\to \widehat u^{-}
       &&\text{in }L^p(\Omega;C(\overline{I}\times\overline{Y^{-}})),\label{eq:inner_strong}\\
 h_\e^{\rm unf}&\to h
       &&\text{in }L^p(\Omega;C(\overline{I};C^{2+\beta'}(\Gamma))),\label{eq:height_strong}
\end{align}
along a subsequence.
Moreover,
\begin{equation}\label{eq:height_velocity}
 \partial_th_\e^{\rm unf}\to\partial_th
       \quad\text{in }L^2(I\times\Omega\times\Gamma)
       \qquad 
       \nabla_yU_\e^-\rightharpoonup\nabla_y\widehat u^-
       \quad\text{in }L^2(I\times\Omega\times Y^-).
\end{equation}
For almost every $x$, the limits are the unique admissible solution of \systemref{sys:driven_cell}{\widehat P_{\rm cell}} with inputs
\[
 (a,g,F)=(u_0^-(x,\cdot),u^{+}(\cdot,x),f^{-}(\cdot,x,\cdot)).
\]
\end{proposition}
\begin{proof}
For $x\in\Omega$, let $k_\e(x)\in\indexI$ denote the index of the cell
containing $x$.\footnote{Since $\dist(x,\partial\Omega)>0$, we have $x\in\e(k_\e(x)+\overline Y)\subset\Omega$ for all sufficiently small $\e$;
of course, the cell index $k_\e(x)$ changes with $\e$ (and is chosen arbitrarily if $x$ lies on a cell boundary).}
We set
\[
 a_\e(x)=\mathcal T_\e u_{\e,0}^-(x,\cdot),\qquad
 g_\e(x)=U_\e^+(\cdot,x,\cdot),\qquad
 F_\e(x)=\mathcal T_\e f_\e^-(\cdot,x,\cdot).
\]
Now, Assumption~\ref{ass:macro_initial}, the uniform cellwise bounds in \ref{ass:initial_data}, and interpolation guarantee $a_\e\to u_0^-$ in $L^2(\Omega;C(\overline{Y^-}))$.
Similarly, $g_\e\to u^+$ in $L^2(\Omega;L^2(I\times\Gamma))$ by  \eqref{eq:outer_trace_strong} and $F_\e\to f^-$ in $L^2(\Omega;L^2(I\times Y))$ by \eqref{eq:source_unfold_limit}.
Passing to a subsequence, these convergences hold for almost every $x\in\Omega$.
We fix such an $x$ and write $k=k_\e(x)$ with $\e$ sufficiently small.
The pair $(U_\e^-(\cdot,x,\cdot),h_\e^{\rm unf}(\cdot,x,\cdot))=(\widehat\mu_{\e,k}^-,h_{\e,k})$ is an admissible solution of \systemref{sys:driven_cell}{\widehat P_{\rm cell}} for the inputs $(a_\e(x),g_\e(x),F_\e(x))$.

By \cref{lem:cell_paths}, these pairs are relatively compact in $C(\overline I\times\overline{Y^-})\times C(\overline I;C^{2+\beta'}(\Gamma))$ with bounded
gradients.
By \cref{lem:cell_closed}, every cluster point is an admissible solution for the inputs
$(u_0^-(x,\cdot),u^+(\cdot,x),f^-(\cdot,x,\cdot))$, and by \cref{thm:cell_stability} this
solution is unique.
Denoting it by $(\widehat u^-(x),h(x))$, the whole sequence converges, and
$\partial_th_\e^{\rm unf}(x)\to\partial_th(x)$ in $L^2(I\times\Gamma)$ by \cref{lem:cell_closed}.

As an almost everywhere limit of measurable functions, $x\mapsto(\widehat u^-(x),h(x))$ is also measurable.
Since all quantities are uniformly bounded, dominated convergence yields \eqref{eq:inner_strong}, \eqref{eq:height_strong}, and \eqref{eq:height_velocity}.
Finally, $\nabla_yU_\e^-$ is bounded in $L^2(I\times\Omega\times Y^-)$ by
\cref{lem:cell_paths} and the unfolding identities, and every weak limit equals
$\nabla_y\widehat u^-$ by \eqref{eq:inner_strong}.
\end{proof}

The strong convergence in \cref{thm:micro_compactness} implies strong convergence of the Hanzawa Jacobians, elliptic matrices, and surface Jacobians in every finite $L^p$ norm; 
since they are also uniformly bounded in $L^\infty$, their products with fixed $L^2$ functions converge strongly in $L^2$ by dominated convergence.
The transformed cell velocities converge strongly in $L^2$, and so do the traces $U_\e^\pm|_\Gamma$, by the trace interpolation inequality and the outer trace convergence, respectively.
In particular,
\[
 J_{h_\e^{\rm unf}}U_\e^{-}\to J_h\widehat u^{-},
 \qquad A_{h_\e^{\rm unf}}\nabla_yU_\e^{-}\rightharpoonup
              A_h\nabla_y\widehat u^{-}
\]
in $L^2(I\times\Omega\times Y^-)$, and
$j_{h_\e^{\rm unf}}(U_\e^--U_\e^+)\to j_h(\widehat u^--\widehat u^+)$
in $L^2(I\times\Omega\times\Gamma)$.

The phase volumes \eqref{eq:phase_volumes} satisfy $\theta^++\theta^-=1$ and they are bounded away from zero with $\theta^\pm(0,\cdot)=|Y^\pm|$.
Moreover, they are uniformly Lipschitz in time with values in $L^\infty(\Omega)$, and
\begin{equation}\label{eq:volume_velocity}
 \partial_t\theta^{-}=\int_\Gamma j_h\beta(h)\,\partial_th\di S_y
 =-\partial_t\theta^{+}.
\end{equation}
Indeed, $\theta^-(t,x)=|Y^-[h(t,x)]|$ and, by the transport theorem, $\partial_t\theta^-$ is the surface integral of the normal velocity
over $\Gamma[h]$ (see, e.g., \cite[Chapter~2]{PruessSimonett2016}).

\subsection{Effective flux and proof of the homogenization theorem}
\label{ssec:hom_flux}
It remains to determine the limiting flux and pass to the macroscopic conservation law.

\begin{lemma}
\label{lem:outer_corrector}
After a further extraction there is $u_1^{+}\in L^2(I\times\Omega;\mathcal V)$ such that
\begin{equation}\label{eq:outer_gradient_limit}
 \mathcal T_\e(\nabla_x\widehat u_\e^{+})
 \rightharpoonup \nabla_xu^{+}+\nabla_yu_1^{+}
       \quad\hbox{in }L^2(I\times\Omega\times Y^{+}).
\end{equation}
For almost every $(t,x)$, it satisfies
\begin{equation}\label{eq:outer_cell_weak}
 \int_{Y^{+}}A_{h(t,x)}
       (\nabla_xu^{+}+\nabla_yu_1^{+})\cdot\nabla_y\psi\di y=0
       \qquad\hbox{for all }\psi\in\mathcal V.
\end{equation}
Consequently,
\begin{equation}\label{eq:corrector_representation}
 u_1^{+}(t,x,\cdot)=\sum_{j=1}^3\partial_{x_j}u^{+}(t,x)\,\chi_j[h(t,x)]
 \qquad\text{for a.e. }(t,x)\in I\times\Omega.
\end{equation}
\end{lemma}
\begin{proof}
We apply \cref{lem:unfold_convergence} to $P_\e^+\widehat u_\e^{+}$.
Along a subsequence, its unfolded gradients converge weakly in $L^2(I\times\Omega\times Y)$ to $\nabla_xu^{+}+\nabla_y\widetilde u_1$ with $\widetilde u_1\in L^2(I\times\Omega;H^1_{\#}(Y)/\R)$.
Restriction to $Y^{+}$ then gives \eqref{eq:outer_gradient_limit}.
We take the test function
\[
 \varphi_\e(t,x)=\e\,\varphi(t,x)\psi(x/\e),
 \quad \varphi\in C_c^\infty(I\times\Omega),
 \quad \psi\in C^\infty_{\#}(Y)
\]
for the outer phase in \eqref{eq:weak_form}.
Due to the additional $\e$, the time and source terms are $\mathcal O(\e)$, as is the contribution of $J_\e b_\e\widehat u_\e^{+}\cdot\nabla \varphi_\e$ because of $b_\e=\mathcal O(\e)$.
Letting $\e\to0$ and using \eqref{eq:outer_gradient_limit} and the strong convergence $\mathcal T_\e(A_\e)\to A_h$ in $L^p(I\times\Omega;L^\infty(Y))$, $p<\infty$, we obtain
\[
 \int_I\int_\Omega\varphi\int_{Y^+}A_h(\nabla_xu^++\nabla_yu_1^+)\cdot\nabla_y\psi
 \di{y}\di{x}\di{t}=0 .
\]
Density of $C^\infty_\#(Y^+)$ in $\mathcal V$ and localization in $(t,x)$ give \eqref{eq:outer_cell_weak}.

Finally, by \cref{lem:cell_correctors}, \eqref{eq:outer_cell_weak} has for
almost every $(t,x)$ a unique solution in $\mathcal V$, and by linearity this
solution is the right-hand side of \eqref{eq:corrector_representation}, which
is measurable in $(t,x)$ by the Lipschitz dependence of $\chi_j[h]$ on $h$.
\end{proof}

\begin{remark}[Physical cell]\label[remark]{rem:physical_tensor}
$D(h)$ is the usual Neumann homogenized tensor of the physical cell
$Y^{+}[h]$, normalized by $|Y|=1$:
\[
 \xi\cdot D(h)\xi=d^{+}\min_{\psi\in H^1_{\#}(Y^{+}[h])}
   \int_{Y^{+}[h]}|\xi+\nabla_y\psi|^2\di y .
\]
This follows from the variational characterization of $D(h)$ by the change of variables $y\mapsto\Theta_h(y)$, which is the identity near $\partial Y$ and thus preserves periodicity.
\end{remark}

\begin{proposition}[Limit of the outer conservation law]
\label{pro:macro_balance}
The zero-extended reference fluxes satisfy
\begin{equation}\label{eq:macro_flux_limit}
 \chi_{\Omega_\e^+}\big(d^{+}A_\e\nabla\widehat u_\e^{+}
                      +J_\e b_\e\widehat u_\e^{+}\big)
       \rightharpoonup D(h)\nabla_xu^{+}
       \quad\hbox{in }L^2(I\times\Omega),
\end{equation}
and for every $\varphi\in C^\infty(\overline I\times\overline\Omega)$ with
$\varphi(T)=0$,
\begin{multline}\label{eq:macro_weak}
 -\int_0^T\!\int_\Omega\theta^{+} u^{+}\partial_t\varphi\di x\di t
 +\int_0^T\!\int_\Omega D(h)\nabla_xu^{+}\cdot\nabla_x\varphi\di x\di t\\
 =\int_0^T\!\int_\Omega\theta^{+} f^{+}\varphi\di x\di t
 +|Y^{+}|\int_\Omega u_0^+\varphi(0)\di x.
\end{multline}
\end{proposition}
\begin{proof}
Since $\|b_\e\|_\infty\le C\e$, the transport part of the flux tends to zero
in $L^2$. By \cref{lem:outer_corrector}, the unfolded diffusive flux converges
weakly in $L^2(I\times\Omega\times Y^+)$ to
$d^{+}A_h\big(\nabla_xu^{+}+\sum_j\partial_{x_j}u^{+}\nabla_y\chi_j[h]\big)$,
whose cell average is $D(h)\nabla_xu^{+}$ by \eqref{eq:effective_tensor}; since
the boundary strip is negligible, this gives \eqref{eq:macro_flux_limit}.

Next we test the outer equation in \eqref{eq:weak_form} with $\varphi$
restricted to $D_\e$. After unfolding, the mass and source terms converge
strongly to $J_hu^{+}$ and $J_hf^{+}$, by \cref{thm:micro_compactness},
\cref{thm:outer_compactness}, and \ref{ass:sources}; their cell averages are
$\theta^{+}u^{+}$ and $\theta^{+}f^{+}$. For the initial term, $J_\e(0)=1$,
$\chi_{D_\e}\stackrel{*}{\rightharpoonup}|Y^+|$ in $L^\infty(\Omega)$, and
$P_\e^+u_{\e,0}^+\to u_0^+$ in $L^2(\Omega)$, so
$\int_{D_\e}u_{\e,0}^+\varphi(0)\di x\to|Y^+|\int_\Omega u_0^+\varphi(0)\di x$.
Together with \eqref{eq:macro_flux_limit}, this gives \eqref{eq:macro_weak}.
\end{proof}

\begin{proof}[Proof of \cref{thm:homogenization}]
By \cref{thm:prepared_schauder}, \cref{pro:basic_estimates}(i), and \cref{thm:outer_compactness}, there is a subsequence along which the outer concentrations converge as in \cref{thm:outer_compactness}; the data converge by \ref{ass:macro_initial} and \ref{ass:macro_sources}.
By \cref{thm:micro_compactness} and \cref{thm:cell_stability}, for almost every $x$ the limits $(\widehat u^-(\cdot,x),h(\cdot,x))$ form the unique solution of the driven cell problem, which gives the inner equation and the surface law with their initial values.
By \cref{lem:outer_corrector} and \cref{pro:macro_balance}, $u^+$ satisfies the outer equation in the weak form \eqref{eq:macro_weak}.
Hence the limit is a weak solution of \systemref{sys:homogenized-fixed}{\widehat P_0} in the sense of \cref{def:weak_homogenized}.

The bounds in (iii) follow from those of \cref{thm:prepared_schauder}:
$L^\infty$ bounds are preserved under a.e.\ convergence, closed balls of
$C^{2+\beta}(\Gamma)$ are closed under convergence in $C^{2+\beta'}(\Gamma)$
for every $t\in\overline I$, and the velocity bound passes to the weak-$*$
limit of $\partial_th_\e^{\rm unf}$.
\end{proof}

\appendix
\section{Auxiliary results}\label{app:auxiliary}
We collect geometric properties of the Hanzawa transformation, extension and trace estimates, and the unfolding operators.

\subsection{The Hanzawa transformation}\label{app:hanzawa}
In the tubular coordinates $y=\gamma+\zeta n_\Gamma(\gamma)$, write $\mathcal E v(y)=\chi(\zeta)v(\gamma)n_\Gamma(\gamma)$ in $U_\Gamma$ and $\mathcal Ev=0$ elsewhere, so that $\Theta_h=\mathrm{Id}+\mathcal E h$.
A finite partition of unity gives $\|\mathcal Ev\|_{W^{k,s}(Y)}\le C\|v\|_{W^{k,s}(\Gamma)}$ for $k=0,1,2$ and $1\le s\le\infty$.

\begin{lemma}
\label{lem:geometric_estimates}
There are $c,C>0$ such that, for all $h$ with $\|h\|_{C^2(\Gamma)}\le\eta^*$,
\[
 c\le J_h\le C,\qquad c\le\beta(h),\qquad
 c|\xi|^2\le A_h\xi\cdot\xi\le C|\xi|^2,
\]
the matrices $D\Theta_h^{\pm1}$, the functions $J_h^{\pm1}$, and $A_h$ are bounded in $W^{1,\infty}(Y)$ by $C$, and $\|N_h\psi-\nabla\psi\cdot n_\Gamma\|_{L^\infty(\Gamma)}\le C\|\psi\|_{C^1(\overline Y)}\|h\|_{C^1(\Gamma)}$ for $\psi\in C^1(\overline Y)$.
If moreover $\partial_th\in C(\Gamma)$, then $\|b_h\|_{L^\infty(Y)}\le C\|\partial_th\|_{C(\Gamma)}$.
\end{lemma}
\begin{proof}
By \eqref{eq:hanzawa_cell}, $\|D\Theta_h-I\|_{W^{1,\infty}}\le C\|h\|_{C^2(\Gamma)}$, and the choice of $\eta^*$ keeps $D\Theta_h$ uniformly close to the identity; matrix inversion and the determinant then give all bounds, and the formula for $\widehat n_h$ gives $\beta(h)\ge c$ and the estimate for $N_h$.
Finally, $\partial_t\Theta_h=\mathcal E(\partial_th)$; see also \cite{PruessSimonett2016}.
\end{proof}

\begin{lemma}\label{lem:app_geometry_flux}
Let $h$ be admissible with $\partial_th\in C(\Gamma)$. 
Then,
\begin{equation*}
 \partial_tJ_h=\dive(J_hb_h)\quad\text{in }Y,\qquad
 J_hb_h\cdot n_\Gamma=j_h\widehat V_h\quad\text{on }\Gamma .
\end{equation*}
Moreover, for $d>0$, a function $u$ on $Y^+[h]$ or $Y^-[h]$, and $z=u\circ\Theta_h$, satisfy
\begin{equation*}
 (dA_h\nabla z+J_hb_hz)\cdot n_\Gamma
 =j_h\bigl(d(\nabla u\cdot n[h])\circ X_h+z\widehat V_h\bigr)\quad\text{on }\Gamma.
\end{equation*}
In particular, since $j_h>0$, $u$ satisfies $d\nabla u\cdot n[h]+uV=0$ on $\Gamma[h]$ if and only if $(dA_h\nabla z+J_hb_hz)\cdot n_\Gamma=0$ on $\Gamma$, and $\partial_tu-d\Delta u=f$ in $Y^\pm[h]$ becomes $\partial_t(J_hz)-\dive(dA_h\nabla z+J_hb_hz)=J_h(f\circ\Theta_h)$ in $Y^\pm$.
\end{lemma}
These identities follow from the cofactor formula, Jacobi's formula, and the Piola identity; see \cite[Section~1.3]{PruessSimonett2016}.

\begin{lemma}\label{lem:app_coefficient_differences}
Let $h_1,h_2$ satisfy $\|h_j\|_{C^2(\Gamma)}\le\eta^*$ and $\|\partial_th_j\|_{L^\infty(\Gamma)}\le M$, and set $\rho=h_1-h_2$, $J_j=J_{h_j}$, $K_j=dA_{h_j}/J_j$, and $c_j=b_{h_j}-K_j\nabla\log J_j$ for fixed $d>0$.
For $1<s<\infty$ and almost every time,
\begin{align*}
 \|K_1-K_2\|_{L^s(Y)}+\|J_1-J_2\|_{L^s(Y)}&\le C\|\rho\|_{W^{1,s}(\Gamma)},\\
 \|J_1-J_2\|_{L^s(\Gamma)}+\|\beta(h_1)-\beta(h_2)\|_{L^s(\Gamma)}
       &\le C\|\rho\|_{W^{1,s}(\Gamma)},\\
 \|c_1-c_2\|_{L^s(Y)}
       &\le C\bigl(\|\rho\|_{W^{2,s}(\Gamma)}+\|\partial_t\rho\|_{L^s(\Gamma)}\bigr),
\end{align*}
with $C$ depending only on $s$, $M$, $d$, and the reference geometry.
\end{lemma}
\begin{proof}
$D\Theta_h$, $J_h$, and $K_h$ are smooth functions of $(h,\nabla_\Gamma h)$ on the admissible set, so the mean-value theorem and the bound on $\mathcal E$ give the first line, and one spatial derivative gives $\|\nabla\log J_1-\nabla\log J_2\|_{L^s}\le C\|\rho\|_{W^{2,s}}$.
For the velocity coefficient, 
\[
b_{h_1}-b_{h_2}=(D\Theta_{h_1})^{-1}\mathcal E(\partial_t\rho)+((D\Theta_{h_1})^{-1}-(D\Theta_{h_2})^{-1})\mathcal E(\partial_th_2),
\]
so no spatial derivative of the velocity is needed.
Combining these bounds in $c_1-c_2=(b_{h_1}-b_{h_2})-(K_1-K_2)\nabla\log J_1-K_2(\nabla\log J_1-\nabla\log J_2)$ gives the third line.
\end{proof}

\begin{lemma}
\label{lem:jac_time_translates}
For every $\bh\in\mathcal H_M(I)$, we have
\begin{equation}\label{eq:jac_time_modulus}
 \|J_\e(t+s)-J_\e(t)\|_{L^\infty(\Omega)}\le Cs^{1/2},
 \qquad 0\le s\le T,\quad 0\le t\le T-s,
\end{equation}
with $C$ depending only on $M$ and the reference geometry.
\end{lemma}
\begin{proof}
In each cell, $\|h(t+s)-h(t)\|_{C^0}\le Ms$, and interpolation with the $C^2$ bound gives $\|h(t+s)-h(t)\|_{C^1}\le Cs^{1/2}$; $J_h$ depends smoothly on $(h,\nabla_\Gamma h)$, and $J_\e=1$ outside the cells.
\end{proof}

\subsection{Approximation and continuity of prescribed motions}
\begin{lemma}\label{lem:app_motion_approximation}
Fix $\e>0$ and $\bh\in\mathcal H_{M}(I)$ with $\bh(0)=0$.
There are smooth $\bh_n$ with $\bh_n(0)=0$, $\bh_n\to\bh$ in $\mathcal X_\e(I)$, $C^1$ bound at most $\eta^*$, $C^2$ bound $\frac{\eta^*}2+o(1)$, and velocity bound $M+o(1)$.
\end{lemma}
\begin{proof}
Extend $h$ oddly to negative times and by reflection past $T$, mollify in space and with a one-sided time mollifier, subtract the value at $t=0$, and multiply by a factor tending to one from below; since $\indexI$ is finite, the scales can be chosen uniformly in $k$.
\end{proof}

\begin{lemma}\label{lem:app_bulk_continuity}
Fix $\e>0$ and let $\bh_j,\bh\in\mathcal H_{M'}(I)$, $M'>0$, with $\bh_j(0)=\bh(0)=0$ and $\bh_j\to\bh$ in $\mathcal X_\e(I)$.
Then the corresponding solutions of \cref{lem:prescribed_bulk} converge in $C(\overline{I};L^2(\Omega_\e^\pm))\cap L^2(I;H^1(\Omega_\e^\pm))$.
In particular, energy and comparison bounds proved for smooth motions pass to $\bh$ via \cref{lem:app_motion_approximation}.
\end{lemma}
\begin{proof}
With the coefficients of \eqref{eq:prescribed_density}, $w_j=J_jz_j-Jz$ has zero initial value and solves $\partial_tw_j-\dive(K_j\nabla w_j+c_jw_j)=\dive G_j+F_j$ with zero total conormal flux, where $G_j=(K_j-K)\nabla(Jz)+(c_j-c)Jz$ and $F_j=J_j\widehat f_j-J\widehat f$.
Uniform convergence of the coefficients and \ref{ass:sources} give $G_j,F_j\to0$ in $L^2$, and the energy estimate for fixed $\e$ gives $w_j\to0$; pointwise bounds pass to the limit along an a.e.\ convergent subsequence.
\end{proof}

\subsection{Extension and trace estimates}
\begin{lemma}
\label{lem:hom_extension}
There are linear operators $P_\e^\pm\colon H^1(\Omega_\e^\pm)\to H^1(\Omega)$ such that
\begin{gather*}
    \|P_\e^\pm z\|_{L^q(\Omega)}\lesssim \|z\|_{L^q(\Omega_\e^\pm)}\ \ (q=2,\infty),\qquad
    \|\nabla P_\e^+ z\|_{L^2(\Omega)}\lesssim \|\nabla z\|_{L^2(\Omega_\e^+)},\\
    \|\nabla P_\e^- z\|_{L^2(\Omega)}\lesssim\|\nabla z\|_{L^2(\Omega_\e^-)}+\e^{-1}\|z\|_{L^2(\Omega_\e^-)}.
\end{gather*}
\end{lemma}
\begin{proof}
For $P_\e^+$, extend from the connected set $Y^+$ into $Y^-$ by a Sobolev extension that preserves constants and is bounded on $L^\infty$, apply it after subtracting the cell average, and rescale cellwise; since nothing changes near the cell faces, the pieces join in $H^1(\Omega)$.
For $P_\e^-$, extend from $Y^-$ into a neighborhood compactly contained in $Y$ and multiply by a cutoff; the cutoff produces the term $\e^{-1}\|z\|_{L^2}$.
\end{proof}

\begin{lemma}\label{lem:scaled_trace}
If $\max_k\|h_k\|_{C^2(\Gamma)}\le\eta^*$, then for every $a>0$ there is $C_a$ with
\[
 \e\|z\|_{L^2(\Gamma_\e(t))}^2
 \le a\e^2\|\nabla z\|_{L^2(\Omega_\e^{\pm}(t))}^2
       +C_{a}\|z\|_{L^2(\Omega_\e^{\pm}(t))}^2
 \qquad \bigl(z\in H^1(\Omega_\e^{\pm}(t))\bigr).
\]
\end{lemma}
This follows from the trace inequality on $Y^-[h_k]$ and on the collar $Y_\Sigma^+[h_k]$, uniformly by \cref{lem:geometric_estimates}, and rescaling.

\subsection{Periodic unfolding}\label{app:unfolding}
We follow \cite{CioranescuEtAl2012}, with time as an additional parameter.
Let $\Omega_\e^{\rm cell}=\bigcup_{k\in\indexI}\e(k+Y)$ and $B_\e=\Omega\setminus\Omega_\e^{\rm cell}$, so that $|B_\e|\le C\e$, and let $D_\e^G=\bigcup_{k}\e(k+G)$ for $G\in\{Y,Y^+,Y^-\}$.
With the cellwise restriction $T_{\e,k}v(y)=v(\e(k+y))$, the assembly of a family $(v_k)_k$ on $G$ or $\Gamma$ is $\mathcal R_\e(v_k)_k(x,y)=v_k(y)$ for $x\in\e(k+Y)$ and $0$ for $x\in B_\e$, and the bulk and surface unfolding operators are $\mathcal T_\e v=\mathcal R_\e(T_{\e,k}v)_k$ and $\mathcal T_\e^\Gamma w=\mathcal R_\e(T_{\e,k}w)_k$.
For $1\le p<\infty$,
\begin{equation*}
 \|\mathcal T_\e v\|_{L^p(\Omega\times G)}=\|v\|_{L^p(D_\e^G)},\quad
 \|\mathcal T_\e^\Gamma w\|_{L^p(\Omega\times\Gamma)}^p=\e\|w\|_{L^p(\Gamma_\e)}^p,\quad
 \nabla_y\mathcal T_\e v=\e\mathcal T_\e(\nabla_xv),
\end{equation*}
and unfolding commutes with traces on $\Gamma$.
For the moving geometry we first pull back and then unfold, so that $\mathcal T_\e J_\e=\chi_{\Omega_\e^{\rm cell}}J_{h_\e^{\rm unf}}$ and $\mathcal T_\e A_\e=\chi_{\Omega_\e^{\rm cell}}A_{h_\e^{\rm unf}}$ on $\Omega\times Y$, and similarly for the surface Jacobian.

\begin{lemma}\label{lem:unfold_convergence}
If $v_\e\to v$ in $L^p(I\times\Omega)$, $1\le p<\infty$, then $\mathcal T_\e v_\e\to v$ in $L^p(I\times\Omega\times Y)$.
If moreover $(v_\e)$ is bounded in $L^2(I;H^1(\Omega))$, then, along a subsequence, $\mathcal T_\e(\nabla_xv_\e)\rightharpoonup\nabla_xv+\nabla_yv_1$ in $L^2(I\times\Omega\times Y)$ for some $v_1\in L^2(I\times\Omega;H^1_{\#}(Y)/\R)$.
\end{lemma}
See \cite[Section~2 and Theorem~3.12]{CioranescuEtAl2012} and \cite{Allaire1992}.
Contributions of the strip $B_\e$ are estimated with the available integrability: if $(g_\e)$ is bounded in $L^p(I\times\Omega)$ for some $p>1$, then $\|g_\e\|_{L^1(I\times B_\e)}\le(T|B_\e|)^{1-1/p}\|g_\e\|_{L^p}\to0$.

\subsection{Cell correctors and effective tensor}

\begin{lemma}\label{lem:cell_correctors}
For every $h\in C^1(\Gamma)$ with $\|h\|_{C^1(\Gamma)}\le\eta^*$, the corrector problems \eqref{eq:corrector_basis} have unique solutions $\chi_j[h]\in\mathcal V$, $j=1,2,3$. 
There are constants $c,C>0$, depending only on $\Gamma$ and $\eta^*$, such that for all such $h,h'$:
\begin{enumerate}[label=(\roman*)]
 \item $\|\chi_j[h]\|_{H^1(Y^+)}\le C$ and
       $\|\chi_j[h]-\chi_j[h']\|_{H^1(Y^+)}\le C\|h-h'\|_{C^1(\Gamma)}$;
 \item $D(h)$ is symmetric with $D(h)\xi\cdot\xi\ge c|\xi|^2$ for all
       $\xi\in\R^3$;
 \item $|D(h)-D(h')|\le C\|h-h'\|_{C^1(\Gamma)}$.
\end{enumerate}
\end{lemma}
\begin{proof}
By \cref{lem:geometric_estimates}, $A_h$ is uniformly elliptic and bounded for admissible $h$, with constants depending only on $\Gamma$ and $\eta^*$. 
Since $Y^+$ is connected, Poincar\'e's inequality holds on $\mathcal V$, so Lax--Milgram gives unique solutions $\chi_j[h]\in\mathcal V$ of \eqref{eq:corrector_basis} with $\|\chi_j[h]\|_{H^1(Y^+)}\le C$.

Testing \eqref{eq:corrector_basis} for $\chi_j[h]$ with $\chi_i[h]$ gives
\[
 D_{ij}(h)=d^+\int_{Y^+}(e_i+\nabla_y\chi_i[h])\cdot A_h(e_j+\nabla_y\chi_j[h])\di y,
\]
which is symmetric, and $\xi\cdot D(h)\xi=d^+\min_{\psi\in\mathcal V}\int_{Y^+}(\xi+\nabla_y\psi)\cdot A_h(\xi+\nabla_y\psi)\di y$.
By uniform ellipticity, the lower bound reduces to $\inf_{\psi\in\mathcal V}\|\xi+\nabla_y\psi\|_{L^2(Y^+)}^2\ge c|\xi|^2$, where the infimum is independent of $h$. 
If it failed, compactness would give a unit vector $\xi$ and a periodic $\psi$ with $\nabla_y\psi=-\xi$ on the connected set $Y^+$; since $Y^-\Subset Y$, every pair of opposite faces lies in $Y^+$, and periodicity forces $\xi=0$, a contradiction.

For admissible $h,h'$, subtracting the corrector equations gives
\[
 \int_{Y^+}A_h\nabla_y(\chi_j[h]-\chi_j[h'])\cdot\nabla_y\psi\di y
 =-\int_{Y^+}(A_h-A_{h'})(e_j+\nabla_y\chi_j[h'])\cdot\nabla_y\psi\di y .
\]
Testing with $\psi=\chi_j[h]-\chi_j[h']$ and using $\|A_h-A_{h'}\|_{L^\infty(Y)}\le C\|h-h'\|_{C^1(\Gamma)}$ from \cref{lem:app_coefficient_differences} yields (i); inserting this into \eqref{eq:effective_tensor} shows (iii).
\end{proof}

\section{A local Hölder estimate}\label{app:holder}
We recall a local boundary Hölder estimate for parabolic equations with bounded measurable coefficients and inhomogeneous conormal boundary data. 
The result is directly based on \cite[Theorem~1.4]{Choi2016}.

Let $D\subset\mathbb R^3$ be a smooth bounded domain, $I=(0,T)$, and let $\Gamma_N\subset\partial D$ be relatively open. 
We consider the linear parabolic problem\footnote{Boundary conditions on $\partial D\setminus\Gamma_N$ are irrelevant for the local estimates considered here.}
\begin{equation}\label{eq:local_holder_problem}
\begin{aligned}
  \partial_s z-\dive(a\nabla z)
    &=w\cdot\nabla z+F
      &&\text{in }I\times D,\\
  a\nabla z\cdot\nu
    &=g
      &&\text{on }I\times\Gamma_N,\\
  z(0)&=z_0
      &&\text{in }D.
\end{aligned}
\end{equation}
Here $a$ is measurable and satisfies
\[
  \kappa|\xi|^2\le a(s,x)\xi\cdot\xi,
  \qquad |a(s,x)|\le\kappa^{-1},
\]
for some $\kappa>0$, while $w$, $F$, and $g$ are bounded and $z_0\in W^{2,\infty}(D)$.
For $x_0\in\Gamma_N$, $s_0\in(0,T]$, and $0<\varrho<\varrho_0$ such that $\partial D\cap B_{\varrho_0}(x_0)\subset\Gamma_N$, we set, for $\sigma\in\{\varrho,\varrho_0\}$,
\begin{alignat*}{2}
  s_\sigma&:={\rm pos}(s_0-\sigma^2),
  &\qquad D_\sigma&:=D\cap B_\sigma(x_0),\\
  Q_\sigma&:=(s_\sigma,s_0)\times D_\sigma,
  &\qquad
  \Sigma_\sigma&:=(s_\sigma,s_0)
      \times(\partial D\cap B_\sigma(x_0)).
\end{alignat*}
We set
\[
  K_0:=
  \|z_0\|_{L^\infty(D_{\varrho_0})}
  +\|\Delta z_0\|_{L^\infty(D_{\varrho_0})}
  +\|\nabla z_0\cdot\nu\|_{L^\infty(\partial D\cap B_{\varrho_0}(x_0))}.
\]

\begin{lemma}
\label{lem:local_holder}
Let $z$ be a weak solution of \eqref{eq:local_holder_problem} satisfying
\[
  z\in L^\infty(I\times D)
  \cap C(\overline I;L^2(D))
  \cap L^2(I;H^1(D)).
\]
Then there exist $\alpha\in(0,1)$ and $C>0$, depending only on $\kappa$, $\|w\|_\infty$, $\varrho$, $\varrho_0$, and the local boundary
geometry of $D$ near $x_0\in\Gamma_N$, but independent of $s_0$ and $T$, such that
\begin{equation}\label{eq:local_holder}
  \|z\|_{C^{\alpha/2,\alpha}(\overline{Q_\varrho})}
  \le C\bigl(
    \|z\|_{L^\infty(Q_{\varrho_0})}
    +\|F\|_{L^\infty(Q_{\varrho_0})}
    +\|g\|_{L^\infty(\Sigma_{\varrho_0})}
    +K_0
  \bigr).
\end{equation}
\end{lemma}

\begin{proof}
\emph{Step 1: Flattening.}
Since $\partial D\cap B_{\varrho_0}(x_0)\subset\Gamma_N$ and $\varrho<\varrho_0$, the smaller domain $D_\varrho$ can be covered by finitely many boundary and interior patches contained in $D_{\varrho_0}$.
We consider a boundary patch; the interior case is standard.
Let $\Psi$ be a smooth boundary chart mapping the local domain onto a half-ball $B^+$ and the corresponding boundary portion of $\Gamma_N$ onto its flat part.
In the transformed coordinates, $\widehat z=z\circ\Psi^{-1}$ satisfies an equation of the form
\[
  \partial_s\widehat z
  -\dive(\widehat a\nabla\widehat z)
  =\widehat w\cdot\nabla\widehat z+\widehat F,
  \qquad
  \widehat a\nabla\widehat z\cdot\nu=\widehat g
  \quad\text{on the flat boundary}.
\]
Indeed, the change of variables introduces a smooth positive Jacobian factor in front of the time derivative.
Dividing by this factor produces an additional bounded drift term. The transformed diffusion matrix remains uniformly elliptic and bounded, while
\[
  \|\widehat w\|_\infty
  \le C(1+\|w\|_\infty),\qquad
  \|\widehat F\|_\infty\le C\|F\|_\infty,\qquad
  \|\widehat g\|_\infty\le C\|g\|_\infty,
\]
where $C$ depends only on the boundary chart and
the ellipticity bounds.
Furthermore, the half-ball $B^+$ is convex and therefore satisfies Assumption~1.1 of \cite{Choi2016}, as well as the extension and trace properties (1.1) and (1.13).

\medskip
\emph{Step 2: The case $s_0\ge\varrho_0^2$.}
We set
\[
  K:=\|z\|_{L^\infty(Q_{\varrho_0})}
     +\|F\|_{L^\infty(Q_{\varrho_0})}
     +\|g\|_{L^\infty(\Sigma_{\varrho_0})}.
\]
The assertion is immediate if $K=0$, so we assume $K>0$.
We apply the local Hölder estimate underlying \cite[Theorem~1.4]{Choi2016} to $\widehat z/K$. 
In there, the coefficients are
\[
  A^{ij}=\widehat a_{ij},\qquad
  A^i=C=F^i=0,\qquad
  B^i=-\widehat w_i,\qquad
  f=\widehat F,
\]
with conormal boundary datum $\widehat g$.
The transformed coefficients satisfy the required ellipticity and boundedness conditions.
Choosing $p_i=q_i=12$ for $i=0,1,2,3$, we have $p_{\min}=q_{\min}=12$ as well as
\[
  \frac{3}{p_{\min}}+\frac{2}{q_{\min}}
  =\frac{5}{12}<1.
\]
Thus the integrability conditions in the paper are satisfied, and we may choose $\alpha_{\mathrm C}=1/2$ in \cite[Theorem~1.4]{Choi2016}.
The normalization gives
\[
  \|\widehat z/K\|_\infty\le 1,\qquad
  \|\widehat F/K\|_\infty
  +\|\widehat g/K\|_\infty\le C.
\]
Moreover, $\widehat w$ is uniformly bounded.
When the theorem is applied on cylinders of fixed size, the mixed-norm quantities entering its constants are therefore controlled independently
of $z$, $F$, $g$, $s_0$, and $T$.
In particular, the admissible height $R_1$, the Hölder exponent $\alpha_0$, and the constant in estimate (1.16) can be chosen uniformly.

Although \cite[Theorem~1.4]{Choi2016} is stated for a full cylinder with conormal conditions on its entire lateral boundary, the argument is local. Indeed, its proof uses spatial cut-off functions supported in the respective cylinders, as also observed for the boundedness estimate in \cite[Remark~1.2]{Choi2016}.
Consequently, the same argument applies on boundary patches where the conormal condition holds.
Since
\[
  s_\varrho-s_{\varrho_0}
  =\varrho_0^2-\varrho^2>0,
\]
the smaller cylinder is separated from the lower time boundary of $Q_{\varrho_0}$.
We choose a fixed radius
\[
  0<r\le
  \min\left\{
    R_1,\,
    c(\varrho_0-\varrho),\,
    \frac12\sqrt{\varrho_0^2-\varrho^2}
  \right\},
\]
where $c>0$ depends only on the boundary charts.
This ensures that all local cylinders used to cover $\overline{Q_\varrho}$ remain within the larger region where the equation and the conormal condition hold.

By estimate (1.16), the Hölder seminorm of $\widehat z/K$ on each smaller cylinder is bounded by a uniform constant with exponent $\alpha_0$.
Combining these local estimates with $\|\widehat z/K\|_\infty\le1$ yields a bound on the entire $\overline{Q_\varrho}$.
Indeed, pairs of sufficiently close points are controlled by a common local cylinder, while the $L^\infty$ bound controls pairs at larger parabolic distance.
Returning to the original coordinates and multiplying by $K$, we obtain \eqref{eq:local_holder} with $\alpha=\alpha_0$.

\medskip
\emph{Step 3: The case $s_0<\varrho_0^2$.}
We extend the solution and coefficients to negative times by setting
\begin{gather*}
  \widetilde z(s):=z_0,\qquad
  \widetilde a(s):=I,\qquad
  \widetilde w(s):=0,\\
  \widetilde F(s):=-\Delta z_0,\qquad
  \widetilde g(s):=\nabla z_0\cdot\nu
  \qquad (s<0).
\end{gather*}
By the definition of $K_0$, the extended data satisfy
\[
  \|\widetilde z\|_{L^\infty(\widetilde Q_{\varrho_0})}
  +\|\widetilde F\|_{L^\infty(\widetilde Q_{\varrho_0})}
  +\|\widetilde g\|_{L^\infty(\widetilde\Sigma_{\varrho_0})}
  \le C\bigl(K+K_0\bigr),
\]
where $\widetilde Q_{\varrho_0}$ and $\widetilde\Sigma_{\varrho_0}$ denote the corresponding cylinders with time interval $(s_0-\varrho_0^2,s_0)$.

The extended function $\widetilde z$ is a weak solution of the same parabolic problem on this larger time interval.
Indeed, for $s<0$ the equation and conormal boundary condition follow directly from the definition of $\widetilde F$ and $\widetilde g$.
Moreover, $z(0)=z_0$ in $L^2(D)$ ensures continuity across $s=0$, so that no additional distributional term appears in the time derivative.
The possible jump of $\widetilde a$ at $s=0$ is admissible, since only measurability of the leading coefficients is required.

After shifting the time variable by $\varrho_0^2$, the extended problem falls under Step~2. 
Applying the resulting estimate and restricting to nonnegative times yields \eqref{eq:local_holder}, with the additional contribution $K_0$.
\end{proof}

\bibliographystyle{siamplain}
\bibliography{references}
\end{document}